\pdfoutput=1

\documentclass{siamltex}

\usepackage{amsfonts,amsmath,amssymb,color,verbatim}
\usepackage{stmaryrd}
\usepackage[mathscr]{eucal}
\usepackage[ruled,vlined]{algorithm2e}
\usepackage{accents}
\usepackage{pgfplots}
\usepackage{marginnote}
\usepackage{float} %Floats
\usepackage{placeins} %Floats!!

\usepackage{wasysym}

\def\R          {\ensuremath{\mathbb R}}

\def\M          {\ensuremath{\mathcal M}}

\def\calV         {\ensuremath{\mathcal V}}

\def\bfzero {\ensuremath{\boldsymbol{\mathrm{0}}}}

\def\vecx         {\ensuremath{\boldsymbol{\mathrm{x}}}}
\def\vecy         {\ensuremath{\boldsymbol{\mathrm{y}}}}

\def\matA          {\ensuremath{\boldsymbol{\mathrm{A}}}}
  
\def\matB          {\ensuremath{\boldsymbol{\mathrm{B}}}}
\def\matC          {\ensuremath{\boldsymbol{\mathrm{C}}}}

\def\matI          {\ensuremath{\boldsymbol{\mathrm{I}}}}
\def\matK          {\ensuremath{\boldsymbol{\mathrm{K}}}}

\def\matP          {\ensuremath{\boldsymbol{\mathrm{P}}}}
\def\matQ          {\ensuremath{\boldsymbol{\mathrm{Q}}}}
\def\matR          {\ensuremath{\boldsymbol{\mathrm{R}}}}
\def\matS          {\ensuremath{\boldsymbol{\mathrm{S}}}}

\def\matU          {\ensuremath{\boldsymbol{\mathrm{U}}}}

\def\matX          {\ensuremath{\boldsymbol{\mathrm{X}}}}
\def\matY          {\ensuremath{\boldsymbol{\mathrm{Y}}}}
\def\matZ          {\ensuremath{\boldsymbol{\mathrm{Z}}}}

\def\tensorsize {\ensuremath{{\mathbb R}^{n_1\times \cdots \times n_d}}}
\def\ten{\mathrm{Ten}}

\def\th{t_{1/2}}

\def\TTrank {\ensuremath{\boldsymbol{\mathrm{r}}}}
\def\matDl {\ensuremath{\boldsymbol{\Delta}_L}}
\def\matDr {\ensuremath{\boldsymbol{\Delta}_R}}
\def\Dl {\ensuremath{{\Delta_L}}}
\def\Dr {\ensuremath{{\Delta_R}}}

\def\endCL{\color{black}}

\newcommand{\ud}{\mathrm{d}}

\newcommand*{\dt}[1]{%
  \accentset{\mbox{\large\bfseries .}}{#1}}

\renewcommand{\top}{\intercal}

\title{Time integration of tensor trains}

\author{Christian Lubich\footnotemark[2] \and Ivan V.~Oseledets\footnotemark[3]\ \footnotemark[5] \and Bart Vandereycken\footnotemark[4]}

\begin{document}

\maketitle

\renewcommand{\thefootnote}{\fnsymbol{footnote}}

\footnotetext[2]{Mathematisches Institut,
       Universit\"at T\"ubingen,
       Auf der Morgenstelle 10,
    D--72076 T\"ubingen,
       Germany. (lubich@na.uni-tuebingen.de)}

\footnotetext[3]{Skolkovo Institute of Science and Technology,
Novaya St.~100, Skolkovo, Odintsovsky district, 143025 
Moscow Region, Russia (i.oseledets@skolkovotech.ru)}
\footnotetext[5]{Institute of Numerical Mathematics,
Gubkina St. 8, 119333 Moscow, Russia}

\footnotetext[4]{Department of Mathematics, Princeton University,  Fine Hall, Princeton NJ 08544, USA. 
(bartv@math.princeton.edu)}

\renewcommand{\thefootnote}{\arabic{footnote}}

\begin{abstract} 
A robust and efficient time integrator for dynamical tensor approximation in the tensor train or matrix product state format is presented. The method is based on splitting the projector onto the tangent space of the tensor manifold. The algorithm can be used for updating time-dependent tensors in the 
given data-sparse tensor train / matrix product state format and for computing an approximate solution to high-dimensional tensor differential equations within this data-sparse format. The formulation, implementation and theoretical properties of the proposed integrator are studied, and numerical experiments with problems from quantum molecular dynamics and with iterative processes in the tensor train format are included.
\end{abstract}



\pagestyle{myheadings} \thispagestyle{plain}
\markboth{CH. LUBICH, I.V. OSELEDETS AND B. VANDEREYCKEN}{TIME INTEGRATION OF TENSOR TRAINS}

\section{Introduction}

There has been much interest lately in the development of data-sparse tensor formats for high-dimensional problems ranging from quantum mechanics to information retrieval; see, e.g., the monograph \cite{Hackbusch:2012} and  references therein. A very promising tensor format is provided by tensor trains (TT) \cite{Oseledets:2011a,ot-tt-2009}, which are also known as matrix product states (MPS) in the theoretical physics literature 
\cite{schollwock-2011}
%\cite{Schollwoeck:2011,Verstraete:2008}. 

In the present paper we deal with the problem of computing an approximation to a time-dependent large tensor $A(t), t_0 \leq t \leq \overline t$ within the TT/MPS format. This includes the situation where the tensor $A(t)$ is known explicitly but in a less data-sparse format and we require an approximation of lower complexity. Alternatively, the tensor $A(t)$ could also be defined implicitly as the solution of a 
tensor differential equation $\dt A = F(t,A)$, where $\dt{\phantom{A}}$ denotes $d/dt$. Such a situation typically arises from a space discretization of a
high-dimensional evolutionary partial differential equation.

% approximations in the TT/MPS format to given time-dependent tensors in a less data-sparse format are to be computed, as well as situations where the time-dependent tensor to be approximated in the TT/MPS format is determined implicitly as the solution of a
% tensor differential equation, which
% would typically arise from a space discretization of a
% high-dimensional evolutionary partial differential equation.

In both situations, such an approximation can be obtained by the principle of dynamical low-rank: Given an approximation manifold $\M$, the desired time-dependent approximation $Y(t) \in \M$ is computed as
\[
 \| \dt Y(t) - \dt A(t) \| = \min \qquad \text{or} \qquad \| \dt Y(t) - F(t,Y(t)) \| = \min,
\]
where $\dt A$ and $F$ are given. This is known as the Dirac--Frenkel time-dependent variational principle in physics; see \cite{Kramer:1981,Lubich:2008}. In our case, $\M$ consists of TT/MPS tensors of fixed rank and its manifold structure and tangent space were studied in \cite{Holtz:2012a}. 
For $\| \cdot \|$ the Euclidean norm, the minimizations from above lead to the following differential equations on $\M$:
\begin{equation}\label{eq:ode1}
\dt Y(t) = P_{Y(t)}\ \dt A(t) \qquad \text{and} \qquad \dt Y(t) = P_{Y(t)}\ F(t,Y(t))
\end{equation}
where $P_{Y}$ is the orthogonal projection onto the tangent space of $\M$ at $Y$ (see \S\ref{sec:orth_proj} for a definition of $P_Y$).  This time-dependent variational principle on fixed-rank TT/MPS manifolds is studied in \cite{Lubich:2013}, where the explicit differential equations are derived and their approximation properties are analyzed. We further refer to \cite{Haegeman:2013} for a discussion of time-dependent matrix product state approximations in the physical literature.

A conceptually related, but technically simpler situation arises in the dynamical low-rank approximation of matrices \cite{Koch:2007}. There, the time-dependent variational principle is applied on manifolds of matrices of a fixed rank, in order to update low-rank approximations to time-dependent large data matrices or to approximate solutions to matrix differential equations by low-rank matrices. The arising differential equations for the low-rank factorization need to be solved numerically, which becomes a challenge in the (often occurring) presence of small singular values in the approximation. While standard numerical integrators such as explicit or implicit Runge--Kutta methods then perform poorly, a novel splitting integrator proposed and studied in \cite{Lubich:2014} shows robustness properties under ill-conditioning that are not shared by any standard numerical integrator. The integrator of \cite{Lubich:2014} is based on splitting the orthogonal projector onto the tangent space of the low-rank matrix manifold. It provides a simple, computationally efficient  update of the low-rank factorization in every time step.

In the present paper we extend the projector-splitting integrator of \cite{Lubich:2014} from the matrix case to the TT/MPS case in the time-dependent approximation \eqref{eq:ode1}.

After collecting the necessary prerequisites on tensor trains / matrix product states in \S2, we study the orthogonal projection onto the tangent space of the fixed-rank TT/MPS manifold in \S3. We show that the projector admits an additive decomposition of a simple structure. In \S4 we formulate the algorithm for the splitting integrator based on the decomposition of the projector. In \S5 we show that this integrator inherits from the matrix case an exactness property that gives an indication of the remarkable robustness of the integrator in the presence of small singular values. In \S6 we discuss details of the implementation and present numerical experiments from quantum dynamics and from the application of the integrator to iterative processes in the TT/MPS format.

% 
% {\bf Still todo}
% 
% Why not do projected integrator in ambient space and project back (or retract). Apart from maybe less efficient (although we should not exaggerate this), this should be less robust to overapproximation: using the setting in the overapproximation section, the projector is not Lipschitz uniformly in $\eps$ (curvature is like $\eps^{-1}$). Hence no stable error propagation and also local error depends of $\eps$. The KLS type integrator is shown to be robust to this $\eps$. In addition, when we integrate for $A(h) - A(0)$ with basic Euler step, this integration is on a flat space (no projector) so this is also robust. Problem? I cannot make a numerical experiment where I see this difference (projected is virtually the same!).
% 
% 
%  
% 
%??? problem description, references, outline of the paper ???
%
%??? Should we formulate directly for complex matrices ; makes connection for MPS easier ???
%
%Point out problem with inverses using lifts

\section{Tensor trains / matrix product states: prerequisites}\label{sec:ode}

We present the tensor train  or matrix product state  formats, together with their normalized representations that we will use throughout the paper. Although our presentation is self-contained, its content is not original and can be found in, e.g, \cite{Oseledets:2011a,Holtz:2012a}. 

\subsection{Notation and unfoldings}\label{subsec:unfolding}\hfill

{\it Norm and inner product of tensors.\/}
The norm of a tensor $X\in\tensorsize$, as considered here, is the Euclidean norm of the vector $\vecx$ that carries the
entries $X(\ell_1,\dots,\ell_d)$ of $X$.
The inner product
$\langle X,Y\rangle$ of two tensors
$X,Y\in\tensorsize$ is the Euclidean inner product of the two corresponding vectors $\vecx$ and $\vecy$.

{\it Unfolding and reconstruction.\/} The $i$th unfolding of a
tensor $X\in\tensorsize$ is the matrix
$\matX^{\langle i \rangle}\in\R^{(n_1\cdots n_i)\times (n_{i+1}\cdots n_d)}$
that aligns all entries $X(\ell_1,\dots,\ell_d)$ with fixed $\ell_1,\dots,\ell_i$ in a row of $\matX^{\langle i \rangle}$, and rows and columns are ordered colexicographically. The inverse of unfolding is reconstructing, which we denote as
$$
X= \ten_i (\matX^{\langle i \rangle}),
$$
that is, the tensor $X\in\tensorsize$ has the $i$th unfolding $\matX^{\langle i \rangle}\in\R^{(n_1\dots n_i)\times (n_{i+1}\dots n_d)}$.
 
{\it TT/MPS format.\/}
A tensor $X\in\tensorsize$ is in the TT/MPS format if there exist \emph{core tensors} $C_i \in \R^{r_{i-1}\times n_i \times r_i}$ with  $r_0=r_d=1$ such that 
\[
X(\ell_1,\dots,\ell_d)  = \sum_{j_1=1}^{r_1} \cdots \sum_{j_{d-1}=1}^{r_{d-1}} C_1(1,\ell_1,j_1) \cdot C_2(j_1, \ell_2, j_2) \cdots  C_d(j_{d-1},\ell_d,1)
\]
for $\ell_i=1,\dots,n_i$ and $i=1,\dots,d$. Equivalently, we have
$$
X(\ell_1,\dots,\ell_d) = \matC_1(\ell_1) \cdots \matC_d(\ell_d),
$$
where the $r_{i-1} \times r_i$ matrices $\matC_i(\ell_i)$ are defined as the slices $C_i(:,\ell_i,:)$.

Observe that $X$ can be parametrized by $\sum_{i=1}^d n_i r_{i-1} r_i \leq dNR^2$ degrees of freedom, where $N = \max\{n_i\}$ and $R = \max\{ r_i \}$. In high-dimensional applications where TT/MPS tensors are practically relevant, $R$ is constant or only mildly dependent on $d$. Hence for large $d$, one obtains a considerable reduction in the degrees of freedom compared to a general tensor of size $N^d$.

{\it Left and right unfoldings.\/} For any core tensor $C_i \in \R^{r_{i-1}\times n_i \times r_i}$, we denote
$$%\begin{align*}
\matC_i^<  = \begin{bmatrix} C_i(:,1,:) \\ \vdots \\ C_i(:,n_i,:)
\end{bmatrix} \in \R^{(r_{i-1}n_i)\times r_i}, \qquad
\matC_i^> = \begin{bmatrix} C_i(:,:,1)^\top \\ \vdots \\ C_i(:,:,r_{i})^\top
\end{bmatrix} \in \R^{(r_{i}n_i)\times r_{i-1}}.
$$%\end{align*}
% Matrix $\matC_i^<$ is called the left unfolding of $C_i$ and denoted $L(C_i)$ in \cite{Holtz:2012a} or the unfolding $C_i^{(1,2)}$ in the first and second modes as defined in \cite{Kolda:2009}. Likewise, $\matC_i^>$ is the right unfolding, denoted as $R(C_i)$ or $C_i^{(3,2)}$.
The matrix $\matC_i^<$ is called the \emph{left unfolding} of $C_i$ and $\matC_i^>$ is the \emph{right unfolding}.

{\it TT/MPS rank.\/} We call a vector $\TTrank = (1, r_1, \ldots, r_{d-1}, 1)$ the \emph{TT/MPS rank} of a tensor $X \in \tensorsize$ if
\[
 \rank \matX^{\langle i \rangle} = r_i, \qquad (i=1,\ldots, d-1).
\]
In case $r_i \leq \min\{ \prod_{j=1}^i n_j, \prod_{j=i+1}^d n_j \}$, this implies that $X$ can be represented in the TT/MPS format with core tensors $C_i \in \R^{r_{i-1}\times n_i \times r_i}$ of full multi-linear rank, that is,
\[
 \rank \matC_i^<  =  r_i \qquad \text{and} \qquad \rank \matC_{i}^> = r_{i-1}, \qquad (i=1,\ldots, d).
\]
In addition, it is known (see \cite[Lem.~4]{Holtz:2012a}) that  for fixed $\TTrank$ such a full-rank condition on the core tensors implies that the set
\begin{align}\label{eq:manifold}
 \M &= \{ X \in \tensorsize \colon \text{TT/MPS rank of $X$ is $\TTrank$} \}
\end{align}
is a smooth embedded submanifold in $\tensorsize$.

% 
% We call an TT/MPS tensor $X$ with rank $\TTrank$ a \emph{full rank}  when
% \[
%  \rank \matX^{\langle i \rangle} =  \rank \matC_i^<  = \rank \matC_{i+1}^> = r_i, \qquad (i=1,\ldots, d)
% \]
% In that case, we also have
% 
% Consider all tensors that belong to the set
% It is known (see \cite{Holtz:2012a}) that $M$ is a smooth embedded submanifold in $\tensorsize$ and each $X \in M$ can be parametrized with core tensors of full multilinear rank
% \[
%  \rank \matX^{\langle i \rangle} =  \rank \matC_i^<  = \rank \matC_{i+1}^> = r_i, \qquad (i=1,\ldots, d).
% \]
% 

% A tensor $X\in\tensorsize$ with entries given in the form
% with coefficient matrices $\matC_i(\ell_i)\in \R^{r_{i-1}\times r_i}$ of full rank, and with
%  $r_0=r_d=1$.  called the \emph{$i$th core} of $X$. 

\begin{figure}
     \centering
          \vspace{0.2cm}
     \def\svgwidth{\columnwidth}
     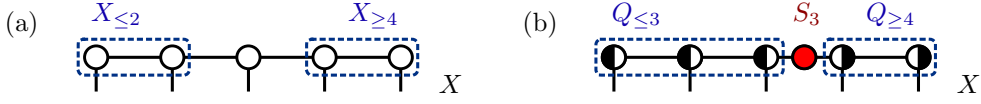
    \caption{A 5 dimensional TT/MPS tensor $X$. Panel (a) indicates specific left and right partial products of $X$. Panel (b) depicts the third recursive SVD of $X$. Observe that left and right orthogonalized cores are denoted using $\LEFTcircle$ and $\RIGHTcircle$ respectively.}\label{fig:TT_left_right}
\end{figure}

{\it Partial products.} Define the left partial product $X_{\le i} \in \R^{n_1 \times \cdots \times n_i \times r_i}$ as
\[
X_{\le i} (\ell_1,\dots,\ell_i,:) = \matC_1(\ell_1) \cdots \matC_i(\ell_i)
\]
and the right partial product $X_{\ge i+1} \in \R^{r_i\times n_{i+1} \times \cdots \times n_d}$ as
\[
X_{\ge i+1} (:,\ell_{i+1},\dots,\ell_d) = \matC_{i+1}(\ell_{i+1}) \cdots \matC_d(\ell_d). 
\]
See also Fig.~\ref{fig:TT_left_right}(a) for their graphical representation in terms of a tensor network.

Let a particular unfolding of each of these partial products be denoted as
$$
\matX_{\le i} \in \R^{(n_1\cdots n_i)\times r_i}, \quad
\matX_{\ge i+1} \in  \R^{(n_{i+1}\cdots n_d)\times r_i}.
$$
The elementwise relation 
$X(\ell_1,\dots,\ell_d)=X_{\le i}(\ell_1,\dots,\ell_i,:) X_{\ge i+1}(:,\ell_{i+1},\dots,\ell_d)$ then translates into
$$
\matX^{\langle i \rangle} = \matX_{\le i} \, \matX_{\ge i+1}^\top.
$$

{\it Recursive construction.\/}  We note the recurrence relations 
\begin{equation}\label{eq:recursive_left}
\matX_{\le i} = ( \matI_{n_i} \otimes \matX_{\le i-1} ) \matC_i^<
\quad\ (i=1,\dots,d)
\end{equation}
starting from $\matX_{\le 0} = 1$, and
\begin{equation}\label{eq:recursive_right}
\matX_{\ge i} = (  \matX_{\ge i+1} \otimes \matI_{n_i} ) \matC_i^>
\quad\ (i=1,\dots,d)
\end{equation}
with $\matX_{\ge d+1}= 1$.  Here $\otimes$ denotes the standard Kronecker product.

Combining the above formulas we note
\begin{equation}\label{eq:matX_i}
\matX^{\langle i \rangle} = (\matI_{n_i} \otimes \matX_{\le i-1}) \matC_i^< \matX_{\ge i+1}^\top,
\end{equation}
which will be an important formula later. Using the recurrence relations for $\matX_{\ge i}$ we also obtain
\begin{equation}\label{eq:matX_imin}
\matX^{\langle i-1 \rangle} = \matX_{\le i-1} \matC_i^{>\top}  (\matX_{\ge i+1} \otimes \matI_{n_i} )^\top,
\end{equation}
which together with the previous formula allows us to pass from the $(i-1)$th to the $i$th unfolding.

\subsection{Left and right orthogonalizations}

Thanks to the recursive relations \eqref{eq:recursive_left} and \eqref{eq:recursive_right}, it is possible to compute the QR decompositions of the matrices $\matX_{\le i}$ and $\matX_{\ge i}$ efficiently. 

Let us explain the case for $\matX_{\le i}$ in detail. First, compute a QR factorization (the ${}^<$ in $\matQ_1^<$ is just notational for now but will become clear in \S\ref{sec:rec_SVD}),
$$
\matX_{\le 1} = \matC_1^< = \matQ_1^< \matR_1, 
\quad\hbox{ with }\quad \matQ_1^{<\top} \matQ_1^< = \matI_{r_1}, \
\matQ_1^< \in \R^{n_1\times r_1},\  \matR_1 \in \R^{r_1\times r_1},
$$
and insert it into the recurrence relation~\eqref{eq:recursive_left} to obtain
$$
\matX_{\le 2} = (\matI_{n_2} \otimes \matQ_1^< \matR_1) \matC_2 ^<
= (\matI_{n_2} \otimes \matQ_1^< )(\matI_{n_2} \otimes  \matR_1)  \matC_2^< .
$$
Next, make another QR decomposition
$$
(\matI_{n_2} \otimes  \matR_1)  \matC_2^< = \matQ_2^< \matR_2, 
\quad\hbox{ with }\quad \matQ_2^{<\top} \matQ_2^< = \matI_{r_2}, \
\matQ_2^<\in \R^{(r_1n_2)\times r_2},\  \matR_2 \in \R^{r_2\times r_2},
$$
so that we have obtained a QR decomposition of
\[
\matX_{\le 2} = \matQ_{\le 2} \matR_2 \qquad  \text{with} \quad \matQ_{\le 2} = (\matI_{n_2} \otimes \matQ_1^< )  \matQ_2^<.
\]
These orthogonalizations can be continued in the same way for $i=2,3,\ldots$. Putting $\matQ_{\le 0} = 1$, we have obtained for each $i=1,\dots,d$ the QR decompositions
\[
 \matX_{\le i}= \matQ_{\le i} \matR_i \qquad \text{with} \qquad \matQ_{\le i}=(\matI_{n_i}\otimes \matQ_{\le i-1}) \matQ_i^<
\]
where the matrices $\matQ_i^< \in \R^{(r_{i-1}n_i)\times r_i}$ and $\matR_i \in \R^{r_i\times r_i}$ are obtained recursively from QR decompositions of lower-dimensional matrices $(\matI_{n_i} \otimes  \matR_{i-1})  \matC_i^< = \matQ_i^< \matR_i$. We call the left partial product $\matX_{\le i}$ in that case \emph{left-orthogonalized}.

% 
% We then obtain
% $
% \matX_{\le 3} = ( \matI_{n_3} \otimes (\matI_{n_2} \otimes \matQ_1 ) \matQ_2) \matQ_3 \matR_3
% $ in the same way
% and further, for $i=1,\dots,d$,
% $$
% \matX_{\le i} = (\matI_{n_i}\otimes(\matI_{n_{i-1}}\otimes \dots (\matI_{n_2}\otimes
% \matQ_1)\dots \matQ_{i-1})\matQ_i \matR_i,
% $$
% with $\matQ_i^{\top} \matQ_i = \matI_{r_i}, \
% \matQ_i\in \R^{(r_{i-1}n_i)\times r_i},\  \matR_i \in \R^{r_i\times r_i}$ constructed from the QR decomposition $(\matI_{n_i} \otimes  \matR_{i-1})  \matC_i^< = \matQ_i \matR_i$.
% This yields a QR decomposition  $\matX_{\le i}= \matQ_{\le i} \matR_i$ in terms of the lower-dimensional matrices $\matQ_1,\dots,\matQ_i$, with $\matQ_{\le i}$ determined recursively as 
% $
% \matQ_{\le i}=(\matI_{n_i}\otimes \matQ_{\le i-1}) \matQ_i.
% $

In a completely analogous way, we can obtain a \emph{right-orthogonalized} $\matX_{\ge i}$ as follows. Denote $\matQ_{\ge d+1} = 1$. Then, starting with $\matX_{\ge d} = \matC_d^> = \matQ_d^> \matR_d$, we can use~\eqref{eq:recursive_right} to obtain the QR decompositions
\begin{equation}\label{eq:Qright}
 \matX_{\ge i}= \matQ_{\ge i} \matR_i \qquad \text{with} \qquad \matQ_{\ge i}=(\matQ_{\ge i+1} \otimes \matI_{n_i} ) \matQ_i^>,
\end{equation}
where the matrices $\matQ_i^> \in \R^{(r_{i}n_i)\times r_{i-1}}$ and $\matR_i \in \R^{r_{i-1}\times r_{i-1}}$ are recursively obtained from $(\matR_{i+1} \otimes \matI_{n_i} )  \matC_i^> = \matQ_i^> \matR_i$. We remark that these $\matR_i$ are in general different than those obtained while orthogonalizing from the left.

% $$
% \matX_{\ge i} = (\matI_{n_i}\otimes(\matI_{n_{i+1}}\otimes \dots (\matI_{n_{d-1}}\otimes
% \matQ_d^>)\dots \matQ_{i+1}^>)\matQ_i^> \matR_i^>,
% $$
% with $\matQ_i^{>\top} \matQ_i^> = \matI_{r_{i-1}}, \
% \matQ_i^>\in \R^{(r_{i}n_i)\times r_{i-1}},\  \matR_i^> \in \R^{r_{i-1}\times r_{i-1}}$.
% This yields a QR decomposition  $\matX_{\ge i}=\matQ_{\ge i}\matR_i^>$.

Observe that when $\matX_{\le i}$ is left-orthogonalized, then so is $\matX_{\le j}$ for any $j < i$. Since $\matX_{\le d} = \matX^{\langle d \rangle}$, we call $X$ \emph{left orthogonal} if $\matX_{\le d}$ is left-orthogonalized. As is evident from Fig.~\ref{fig:TT_left_right}, such a left orthogonal $X$ is recursively computed by modifying the cores $C_i$ from left to right during a so-called \emph{forward sweep}. Likewise, we call $X$ \emph{right orthogonal} if $\matX_{\ge 1} = \matX^{\langle 1 \rangle} $ is right-orthogonalized which is obtained by a \emph{backward sweep} from right to left.

\subsection{Recursive SVD}\label{sec:rec_SVD}

Suppose that $\matX_{\le i} = \matQ_{\le i}  \matR_i$ and $\matX_{\ge i+1} = \matQ_{\ge i+1} \matR_{i+1}$ are QR decompositions obtained from left and right orthogonalizations, we then have the following SVD-like decomposition
\begin{equation}\label{eq:rec_SVD}
 \matX^{\langle i \rangle} = \matQ_{\le i} \matS_i \matQ_{\ge i+1}^\top, \qquad \text{with} \quad  \matS_i= \matR_i \matR_{i+1}^\top \in \R^{r_i \times r_i}.
\end{equation}
The matrix $\matS_i$ can be chosen diagonal, although we do not insist that it is. Since the orthonormal matrices $\matQ_{\le i}$ and $\matQ_{\ge i+1}$ satisfy the recursive relations as explained before, we call \eqref{eq:rec_SVD} a \emph{recursive SVD of $\matX^{\langle i \rangle}$}, or the $i$th recursive SVD of $X$. The graphical representation of such a recursive SVD is depicted in Fig.~\ref{fig:TT_left_right}(b).

This recursiveness can be used for the SVD of $\matX^{\langle i+1 \rangle}$. By~\eqref{eq:Qright}, we can write
\begin{equation}\label{eq:rec_svd_i}
 \matX^{\langle i \rangle} = (\matQ_{\le i} \matS_i) \matQ_{i+1}^{>\top} (\matQ_{\ge i+2} \otimes \matI_{n_{i+1}})^\top.
\end{equation}
To obtain a decomposition of $\matX^{\langle i+1 \rangle}$ by means of the relations \eqref{eq:matX_imin} and \eqref{eq:matX_i}, we identify \eqref{eq:rec_svd_i} with \eqref{eq:matX_imin} (hence, $i-1$ takes the role of $i$ and $\matC_i^{>\top}$ that of $\matQ_{i+1}^{>\top}$). The corresponding expression for \eqref{eq:matX_i} then becomes
% \begin{equation}\label{eq:rec_svd_i}
%  \matX^{\langle i \rangle} = \matQ_{\le i} \matS_i \matQ_{\ge i+1}^\top = (\matQ_{\le i} \matS_i) \matQ_{i+1}^{>\top} (\matQ_{\ge i+2} \otimes \matI_{n_{i+1}})^\top
% \end{equation}
% also implies the unfolding
\begin{equation}\label{eq:rec_svd_iplusA}
 \matX^{\langle i+1 \rangle} = (\matI_{n_{i+1}} \otimes \matQ_{\le i}\matS_i) \matQ_{i+1}^< \matQ_{\ge i+2}^\top,
\end{equation}
which we can also write as
\begin{equation}\label{eq:rec_svd_iplus}
 \matX^{\langle i+1 \rangle} = (\matI_{n_{i+1}} \otimes \matQ_{\le i}) (\matI_{n_{i+1}} \otimes \matS_i) \matQ_{i+1}^< \matQ_{\ge i+2}^\top .
\end{equation}
% Interpret the factor $\matQ_{i+1}^>$ from the QR decomposition as the right unfolding of a tensor $Q_{i+1} \in \R^{r_{i} \times n_{i+1} \times r_{i+1}}$.  By identifying
% \begin{equation}\label{eq:rec_svd_i}
%  \matX^{\langle i \rangle} = \matQ_{\le i} \matS_i \matQ_{\ge i+1}^\top = (\matQ_{\le i} \matS_i) \matQ_{i+1}^{>\top} (\matQ_{\ge i+2} \otimes \matI_{n_{i+1}})^\top
% \end{equation}
% with \eqref{eq:matX_imin}, we obtain using \eqref{eq:matX_i} and the left unfolding of $Q_{i+1}$ that
% \begin{equation}\label{eq:rec_svd_iplus}
%  \matX^{\langle i+1 \rangle} = (\matI_{n_{i+1}} \otimes \matQ_{\le i}\matS_i) \matQ_{i+1}^< \matQ_{\ge i+2}^\top = (\matI_{n_{i+1}} \otimes \matQ_{\le i}) (\matI_{n_{i+1}} \otimes \matS_i) \matQ_{i+1}^< \matQ_{\ge i+2}^\top .
% \end{equation}
Hence, after a QR decomposition
    \begin{equation}\label{eq:rec_SQ_QS}
     (\matI_{n_{i+1}} \otimes \matS_i) \matQ_{i+1}^< = \overline \matQ_{i+1}^< \overline \matS_{i+1},
    \end{equation}    
we obtain the $(i+1)$th recursive SVD
    \[
     \matX^{\langle i+1 \rangle} = \overline  \matQ_{\le i+1} \overline \matS_{i+1} \matQ_{\ge i+2}^\top \qquad \text{with} \quad  \overline  \matQ_{\le i+1}  = (\matI_{n_{i+1}} \otimes \matQ_{\le i}) \overline \matQ_{i+1}^<.
    \]    
    
A similar relation holds between $\matX^{\langle i \rangle}$ and $\matX^{\langle i-1 \rangle}$. Let
% \begin{equation}\label{eq:rec_svd_i_b}
%  \matX^{\langle i \rangle} =  \matQ_{\le i} \matS_i \matQ_{\ge i+1}^\top =  (\matI_{n_{i}} \otimes \matQ_{\le i-1}) \matQ_i^< (\matQ_{\ge i+1}\matS_i^\top)^\top
% \end{equation}
% to
% \begin{equation}\label{eq:rec_svd_imin}
%  \matX^{\langle i-1 \rangle} = \matQ_{\le i-1} \matQ_i^{>\top} ( \matS_i \otimes \matI_{n_{i}}) (\matQ_{\ge i+1}^\top \otimes \matI_{n_{i}}).
% \end{equation}
\begin{equation}\label{eq:rec_svd_i_b}
 \matX^{\langle i \rangle} =  \matQ_{\le i} \matS_i \matQ_{\ge i+1}^\top  =  (\matI_{n_{i}} \otimes \matQ_{\le i-1}) \matQ_i^< (\matQ_{\ge i+1}\matS_i^\top)^\top,
\end{equation}
then using the QR decomposition
    \begin{equation*}
     (\matS_i^\top  \otimes \matI_{n_{i}}) \matQ_{i}^> = \overline \matQ_{i}^> \overline \matS_{i}^\top
    \end{equation*}    
we can write
    \begin{equation}\label{eq:rec_svd_imin} 
\matX^{\langle i-1 \rangle} = \matQ_{\le i-1} \overline\matS_{i}  \overline\matQ_{\ge i}^\top \qquad \text{where} \qquad   \overline \matQ_{\ge i}  = (\matQ_{\ge i+1}\otimes \matI_{n_{i}} ) \overline \matQ_{i}^>.
    \end{equation}

\section{Orthogonal projection onto the tangent space}\label{sec:orth_proj} Let $\M$ be the embedded manifold of tensors of a given TT/MPS rank $\TTrank$; see \eqref{eq:manifold}. In this section, we derive an explicit formula for the orthogonal projection onto the tangent space $T_X\M$ at $X\in\M$,
$$
P_X:\tensorsize \to T_X\M.
$$
With the Euclidean inner product, the projection $P_X(Z)$ for arbitrary $Z\in\tensorsize$ has the following equivalent variational definition:
$$
\langle P_X(Z), \delta X \rangle =\langle Z, \delta X \rangle \qquad \forall\,\delta X \in T_X\M.
$$

Before we state the theorem, we recall a useful parametrization of $T_X \M$ as introduced in~\cite{Holtz:2012a}. Let $X \in \M$ be left orthogonal, that is, in the decompositions
\begin{equation*}
\matX^{\langle i \rangle} =(\matI_{n_i}\otimes \matX_{\le i-1})\, \matC_i^< \, \matX_{\ge i+1}^\top,
\end{equation*}
the matrices satisfy for all $i=1,\ldots, d-1$
\begin{equation}\label{eq:left_orth_X}
\matX_{\le i}^\top \matX_{\le i} = \matI_{r_i} \quad \text{and} \quad \matC_i^{<\top} \matC_i^< = \matI_{r_i}.
\end{equation}
Define then for $i=1,\ldots,d-1$ the subspaces
$$
\calV_i = \left\{ \ten_i \left[ (\matI_{n_i} \otimes \matX_{\le i-1}) \, \delta\matC_i^< \, \matX_{\ge i+1}^\top  \right] \colon \delta C_i \in \R^{r_{i-1}\times n_i \times r_i} \  \text{and} \ \matC_i^{<\top} \delta\matC_i^< = \bfzero \right\} 
$$
and also the subspace
$$
\calV_d = \left\{ \ten_d \left[ (\matI_{n_d} \otimes \matX_{\le d-1}) \, \delta\matC_d^<  \right] \colon \delta C_d \in \R^{r_{d-1}\times n_d \times r_d}  \right\} .
$$
Observe that these subspaces represent the first-order variations in $C_i$ in all the representations \eqref{eq:left_orth_X} together with the so-called \emph{gauge conditions} $\matC_i^{<\top}\, \delta\matC_i^<=\bfzero$ when $i \ne d$; there is no gauge condition for $i=d$.  Now, \cite[Thm.~4]{Holtz:2012a} states that
\[
 T_X \M = \calV_1 \oplus \calV_2 \oplus \cdots \oplus \calV_d.
\]
In other words, every $\delta X \in T_X \M$ admits the unique orthogonal\footnote{The orthogonality of the $\calV_i$ spaces is only implicitly present in \cite[Thm.~4]{Holtz:2012a}; it is however not difficult to prove it explicitly thanks to the left-orthogonalization and the gauge conditions.} decomposition
\[ 
 \delta X =  \sum_{i=1}^d \delta X_i , \qquad \text{with} \quad \delta X_i \in \calV_i, \quad \text{and} \quad \langle \delta X_i,\, \delta X_j \rangle  = \delta_{ij}.
\]

Now we are ready to state the formula for $P_X$. It uses the orthogonal projections onto the range of $\matX_{\le i}$, denoted as $\matP_{\le i}$, and onto the range of $\matX_{\ge i}$, denoted as $\matP_{\ge i}$. With the QR decompositions $\matX_{\le i}=\matQ_{\le i}\matR_{i}$  and $\matX_{\ge i}=\matQ_{\ge i}\matR'_{i}$ , these projections become
$$
\matP_{\le i} = \matQ_{\le i}\matQ_{\le i}^\top, \qquad\hbox{and}\qquad
\matP_{\ge i} = \matQ_{\ge i}\matQ_{\ge i}^\top.
$$
We set $\matP_{\le 0}=1$ and $\matP_{\ge d+1}=1$.

\begin{theorem} \label{thm:proj}
    Let $\M$ be the manifold  of fixed rank TT/MPS tensors. Then, the orthogonal projection onto the tangent space of $\M$ at 
$X \in \M$ is given by
\begin{align*}
P_X(Z) = &\sum_{i=1}^{d-1}  \ten_i\bigl[ (\matI_{n_i}\otimes \matP_{\le i-1}) \matZ^{\langle i \rangle} \matP_{\ge i+1} -  \matP_{\le i} \matZ^{\langle i \rangle} \matP_{\ge i+1}  \bigr]
\\
& +  \ten_d \bigl[ (\matI_{n_d}\otimes \matP_{\le d-1}) \matZ^{\langle d \rangle} \bigr] 
\end{align*}
for any $Z\in\tensorsize$.
\end{theorem}

\begin{proof} 
    % We start from the formula
% $$
% \matX^{\langle i \rangle} =(\matI_{n_i}\otimes \matX_{\le i-1})\, \matC_i^< \, \matX_{\ge i+1}^\top,
% $$
% where we can suppose, after orthogonalizations from the left and right, that
% $$
% \matX_{\le i-1}^\top\matX_{\le i-1}=\matI_{r_{i-1}}, \quad\ 
% \matX_{\ge i+1}^\top \matX_{\ge i+1} = \matI_{r_i}.
%  %\matC_i^{<\top}  \matC_i^< = \matI_{n_i} \qquad (i=1,\dots, d-1).
% $$
% It is known from \cite[Thm.~2]{Holtz:2012a} that tangent vectors $\delta X\in T_X\M$ have the representation
% $$
% \delta X = \sum_{i=1}^d \delta X_i, \qquad \delta\matX_i^{\langle i \rangle} = (\matI_{n_i} \otimes \matX_{\le i-1}) \, \delta\matC_i^< \, \matX_{\ge i+1}^\top
% $$
% with coefficient matrices $\delta\matC_i^<$ that are uniquely determined under the gauge conditions
% $$
% \matC_i^{<\top}\, \delta\matC_i^< = \bfzero \qquad (i=1,\dots,d-1).
% $$
% There is no gauge condition for $i=d$. 
We assume that $X$ is given by \eqref{eq:left_orth_X}. For given $Z\in\tensorsize$, we aim to determine $\delta U =P_X(Z)\in T_X\M$ such that
\begin{equation}\label{eq:var_eq}
\langle \delta U, \delta X \rangle = \langle Z, \delta X \rangle \qquad
\forall \, \delta X \in T_X\M.
\end{equation}
Writing $\delta U=\sum_{j=1}^d \delta U_j$ with $\delta U_j \in \calV_j$, this means that we need to determine matrices $\delta \matB_j^<$ in the unfoldings
$$
\delta \matU_j^{\langle j \rangle} = (\matI_{n_j} \otimes \matX_{\le j-1}) \, \delta \matB_j^< \, \matX_{\ge j+1}^\top,
$$
such that the gauge conditions are satisfied
$$
\matC_j^{<\top}\, \delta \matB_j^< = \bfzero \qquad (j=1,\dots,d-1).
$$

Fix an $i$ between $1$ and $d$. Since $\calV_i$ is orthogonal to $\calV_j$ when $j\ne i$, choosing any $\delta X = \delta X_i \in \calV_i$ in \eqref{eq:var_eq} implies
\begin{equation}\label{eq:var_eq2}
\langle \delta U_i, \delta X_i \rangle = \langle Z, \delta X_i \rangle \qquad
\forall \, \delta X_i \in \calV_i.
\end{equation}
Parametrize $\delta X_i \in \calV_i$ as 
$$
\delta\matX_i^{\langle i \rangle} = (\matI_{n_i} \otimes \matX_{\le i-1}) \, \delta\matC_i^< \, \matX_{\ge i+1}^\top
$$
with $ \delta\matC_i^<$ satisfying the gauge condition for $i\ne d$. Then, the left-hand side of \eqref{eq:var_eq2} becomes
\begin{align*}
\langle \delta U_i, \delta X_i \rangle &= \langle \delta \matU_i^{\langle i \rangle}, \delta\matX_i^{\langle i \rangle} \rangle
\\
&= \langle (\matI_{n_i} \otimes \matX_{\le i-1}) \, \delta \matB_i^< \, \matX_{\ge i+1}^\top, (\matI_{n_i} \otimes \matX_{\le i-1}) \, \delta\matC_i^< \, \matX_{\ge i+1}^\top \rangle
\\
&= \langle  \delta \matB_i^< \matX_{\ge i+1}^\top \matX_{\ge i+1} , \delta\matC_i^<  \rangle,
\end{align*}
since $X$ is left orthogonal. Likewise, for the right-hand side we get
\begin{align*}
\langle Z, \delta X_i \rangle &= \langle  \matZ^{\langle i \rangle}, \delta\matX_i^{\langle i \rangle} \rangle
\\
&=\langle (\matI_{n_i} \otimes \matX_{\le i-1})^\top \matZ^{\langle i \rangle} \matX_{\ge i+1}, 
\delta\matC_i^< \rangle.
\end{align*}
Hence, for all matrices $\delta\matC_i^<$ satisfying the  gauge conditions, we must have
$$
\langle  \delta \matB_i^< \matX_{\ge i+1}^\top \matX_{\ge i+1} , \delta\matC_i^< \rangle = \langle (\matI_{n_i} \otimes \matX_{\le i-1})^\top \matZ^{\langle i \rangle} \matX_{\ge i+1}, 
\delta\matC_i^< \rangle,
$$
which implies, with $\matP_i^<$ the orthogonal projector onto the range of $\matC_i^<$ for $i=1,\dots,d-1$ and with $\matP_i^<=\bfzero$ for $i=d$,
$$
 \delta \matB_i^< = (\matI_i - \matP_i^<) (\matI_{n_i} \otimes \matX_{\le i-1})^\top \matZ^{\langle i \rangle} \matX_{\ge i+1} (\matX_{\ge i+1}^\top \matX_{\ge i+1})^{-1},
$$
where $\matI_{i} = \matI_{n_ir_{i-1}}$.
Inserting this expression into the formula for $\delta\matU_i^{\langle i \rangle}$ gives us
$$
\delta\matU_i^{\langle i \rangle} = (\matI_{n_i} \otimes \matX_{\le i-1}) \, (\matI_i - \matP_i^<) (\matI_{n_i} \otimes \matX_{\le i-1})^\top \matZ^{\langle i \rangle} \matX_{\ge i+1} \, (\matX_{\ge i+1}^\top \matX_{\ge i+1})^{-1} \,\matX_{\ge i+1}^\top.
$$
Since $\matP_{\le i-1}=\matX_{\le i-1}\matX_{\le i-1}^\top$, $\matP_{\le i}= (\matI_{n_i} \otimes \matX_{\le i-1}) \matP_i^< (\matI_{n_i} \otimes \matX_{\le i-1})^\top$ and $\matP_{\ge i+1}=\matX_{\ge i+1} (\matX_{\ge i+1}^\top \matX_{\ge i+1})^{-1} \matX_{\ge i+1}^\top$,
%Recalling that $(\matI_{n_i} \otimes \matX_{\le i-1})\matC_i^< = \matX_{\le i}$ and
%$\matX_{\le i}\matX_{\le i}^\top=\matP_{\le i}$, 
this simplifies to
\begin{align*}
\delta\matU_i^{\langle i \rangle} &= (\matI_{n_i} \otimes \matP_{\le i-1} - \matP_{\le i}) \matZ^{\langle i \rangle}
\matP_{\ge i+1} \qquad (i=1,\dots,d-1),
\\
\delta\matU_d^{\langle d \rangle} &= (\matI_{n_d} \otimes \matP_{\le d-1} ) \matZ^{\langle i \rangle}.
\end{align*}
Now $\delta U=\sum_{i=1}^d \delta U_i$ satisfies the projection condition \eqref{eq:var_eq}.
\end{proof}

Although the formula in Theorem~\ref{thm:proj} lends itself well to practical implementation, its cumbersome notation is a nuisance. We therefore introduce a simpler notation for the forthcoming derivations. 

\begin{corollary}\label{cor:eq_PX}
    For $i=0,\ldots,d+1$, define the orthogonal projectors
\begin{align*}
 &P_{\le i}\colon \R^{n_1 \times \cdots \times n_d} \to T_X \M, \  Z \mapsto \ten_i (  \matP_{\le i} \matZ^{\langle i \rangle})\\
 &P_{\ge i}\colon \R^{n_1 \times \cdots \times n_d} \to T_X \M, \  Z \mapsto \ten_{i-1} (  \matZ^{\langle i-1 \rangle}  \matP_{\ge i} )  .
\end{align*}
Then, the projector $P_X$ in Theorem~\ref{thm:proj} satisfies   
\[
 P_X =  \sum_{i=1}^{d-1}  (P_{\le i-1} P_{\ge i+1} - P_{\le i} P_{\ge i+1})  + P_{\le d-1} P_{\ge d+1}.
\]
In addition, $P_{\le i}$ and $P_{\ge j}$ commute for $i < j$.
\end{corollary}
\begin{proof}
    The fact that $P_{\le i}$ commutes with $P_{\ge j}$ follows from the observation that for any $Z \in \R^{n_1 \times \cdots \times n_d}$, $P_{\le i}(Z)$ acts on the rows of $\matZ^{\langle i \rangle}$---and hence also on the rows of $\matZ^{\langle j \rangle}$---while $P_{\ge j}(Z)$ acts on the columns of $\matZ^{\langle j \rangle}$.

    To write $P_X$ using the new notation, we need to work out the term 
    \begin{align*}
     P_{\ge i+1} (P_{\le i-1}(Z)) &= P_{\ge i+1} \bigl[ \ten_{i-1}(  \matP_{\le i-1} \matZ^{\langle i-1 \rangle} ) \bigr] = \ten_i \bigl[ \matY^{\langle i \rangle} \matP_{\ge i+1} \bigr],
    \end{align*}
    with $Y = \ten_{i-1}(  \matP_{\le i-1} \matZ^{\langle i-1 \rangle} )$. Denote the mode-1 matricization of a tensor by ${\cdot}^{(1)}$; see \cite[\S2.4]{Kolda:2009} for a definition. Then, define the tensors $\widehat Z$ and $\widehat Y$, both of size $(n_1\cdots n_{i-1}) \times n_i \times (n_{i+1}\cdots n_d)$, such that $\widehat \matZ{}^{(1)} = \matZ^{\langle i-1 \rangle}$ and $\widehat \matY{}^{(1)} = \matY^{\langle i-1 \rangle}$.
In addition, let $\times_1$ denote the mode-1 multilinear product of a tensor with a matrix; see \cite[\S2.5]{Kolda:2009}. Then, using \cite[p.~426]{Kolda:2009} to compute matricizations of multilinear products, we get
\[
 ( \widehat Z \times_1 \matP_{\le i-1} )^{(1)} =  \matP_{\le i-1} \widehat \matZ{}^{(1)} = \matP_{\le i-1} \matZ^{\langle i-1 \rangle} = \matY^{\langle i-1 \rangle} = \widehat \matY{}^{(1)}.
\]
Hence, we see that $\widehat Y = \widehat Z \times_1 \matP_{\le i-1}$. Using the notation $\cdot^{(1,2)} = \cdot^{(3)T}$ (see again \cite[\S2.4]{Kolda:2009}), we obtain
\[
 \widehat \matY{}^{(1,2)}  = ( \widehat Z \times_1 \matP_{\le i-1} )^{(1,2)} =  (\matI_{n_i} \otimes \matP_{\le i-1}) \widehat \matZ{}^{(1,2)}.
\]
Now, observe that because of the colexicographical ordering of unfoldings and matricizations, we have $\widehat \matZ{}^{(1,2)} = \matZ^{\langle i \rangle}$ and $\widehat \matY{}^{(1,2)} = \matY^{\langle i \rangle}$ and this gives
\[
 P_{\ge i+1} (P_{\le i-1}(Z)) = \ten_i \bigl[ \matY^{\langle i \rangle} \matP_{\ge i+1} \bigr] = \ten_i \bigl[ (\matI_{n_i} \otimes \matP_{\le i-1}) \matZ^{\langle i \rangle}\matP_{\ge i+1} \bigr]. 
\]
The term $P_{\le i} P_{\ge i+1}$ is straightforward, and this finishes the proof.
\end{proof}

\section{Projector-splitting integrator}\label{sec:abstract_splitting}
We now consider the main topic of this paper: a numerical integrator for the dynamical TT/MPS approximation 
\begin{equation}\label{ode}
\dt Y(t) = P_{Y(t)}(\dt A(t)), \qquad Y(t_0)=Y_0\in\M
\end{equation}
of a given time-dependent tensor $A(t)\in\tensorsize$. 

Our integrator is a Lie--Trotter splitting of the vector field $P_Y(\dt A)$. The splitting itself is suggested by the sum in Corollary~\ref{cor:eq_PX}: using $Y$ in the role of $X$, we can write 
\[
 P_Y(\dt A) = P_1^+(\dt A) - P_1^-(\dt A) + P_2^+(\dt A) - P_2^-(\dt A) + \cdots  - P_{d-1}^-(\dt A) + P_d^+(\dt A)
\]
with the orthogonal projectors
\begin{align}
 P_i^+(Z) &= P_{\le i-1} \, P_{\ge i+1} (Z) = \ten_i \bigl[ (\matI_{n_i} \otimes \matP_{\le i-1}) \matZ^{\langle i \rangle}\matP_{\ge i+1} \bigr],  &  (1 \leq i \leq d), \label{eq:def_Pi_plus} \\ 
  P_i^-(Z) &=  P_{\le i} \, P_{\ge i+1}(Z) = \ten_i \bigl[ \matP_{\le i} \matZ^{\langle i \rangle}\matP_{\ge i+1} \bigr] ,  &  (1 \leq i\leq d-1) . \label{eq:def_Pi_min}
\end{align}
By standard theory (see, e.g., \cite[II.5]{Hairer:2006}), any splitting of this sum results in a first-order integrator, and composing it with the adjoint gives a second-order integrator, also known as the Strang splitting. Somewhat remarkably, we shall show in Thm.~\ref{thm:sol_split} that these split differential equations can be solved in closed form. Furthermore, if they are solved from left to right (or from right to left), the whole scheme can be implemented very efficiently. % In section, we show that this integrator is also robust numerically.

\subsection{Abstract formulation and closed-form solutions}\label{sec:forward_sweep}
Let $t_1 - t_0 > 0$ be the step size. One full step of the splitting integrator solves in consecutive order the following initial value problems over the time interval $[t_0,t_1]$:
\begin{alignat*}{2}
\dt Y_1^+ &= +P_1^+( \dt A), & \qquad Y_1^+(t_0) &= Y_0; \\    
 \dt Y_1^- &= -P_1^-( \dt A), & Y_1^-(t_0) &= Y_{1}^+(t_1); \\    
& \ \ \vdots \\
 \dt Y_i^+ &= +P_i^+( \dt A), & Y_i^+(t_0) &= Y_{i-1}^-(t_1); \\
 \dt Y_i^- &= - P_i^-( \dt A), & Y_i^-(t_0) &= Y_{i}^+(t_1); \\
 & \ \ \vdots \\
 \dt Y_d^+ &= +P_d^+( \dt A), & Y_d^+(t_0) &= Y_{d-1}^-(t_1).
\end{alignat*}
Here, $Y_0 = Y(t_0)$ is the initial value of \eqref{ode} and $Y_d^+(t_1)$ is the final approximation for $Y(t_1)$. Observe that one full step consists of $2d-1$ \emph{substeps}.

We remark that the projectors $P_i^+, P_i^-$ depend on the current value of $Y_i^+(t)$ or $Y_i^-(t)$; hence, they are in general time-dependent. For notational convenience, we do not denote this dependence explicitly since the following result states we can actually take them to be time-independent as long as they are updated after every substep. In addition, it shows how these substeps  can be solved in closed  form.
\begin{theorem}   \label{thm:sol_split}  
Let $\Delta A = A(t_1) - A(t_0)$. %For $t_1 - t_0 > 0$  sufficiently small, 
The initial value problems from above satisfy
     \[
      Y_i^+(t_1) = Y_i^+(t_0) + P_i^+  (\Delta A) \qquad \text{and} \qquad      Y_i^-(t_1) = Y_i^-(t_0) - P_i^- (\Delta A),
     \]     
where $P_i^+$ and $P_i^-$ are the projectors at $Y_i^+(t_0)$ and  $Y_i^-(t_0)$, respectively. 
 
In particular, if $Y_i^+(t_0)$ has the recursive SVD
\begin{equation*}
 [Y_i^+(t_0)]^{\langle i \rangle}  = \matQ_{\le i} \matS_{i} \matQ_{\ge i+1}^\top = (\matI_{n_i} \otimes \matQ_{\le i-1}) \matQ_{i}^< \matS_{i} \matQ_{\ge i+1}^\top,
\end{equation*}
with $\matQ_{\le 0} =  \matQ_{\ge d+1} = 1$, then
\begin{equation*}
 [Y_i^+(t_1)]^{\langle i \rangle}  = (\matI_{n_i} \otimes \matQ_{\le i-1}) \, \left\{ \matQ_{i}^< \matS_{i} \ + (\matI_{n_i} \otimes \matQ_{\le i-1}^\top) [\Delta A]^{\langle i \rangle}   \matQ_{\ge i+1} \right\} \, \matQ_{\ge i+1}^\top .
\end{equation*}
Likewise, if $Y_i^-(t_0)$ has the recursive SVD
\[
 [Y_i^-(t_0)]^{\langle i \rangle}  = \matQ_{\le i} \matS_{i} \matQ_{\ge i+1}^\top, 
\]
then
\[
 [Y_i^-(t_1)]^{\langle i \rangle}  =  \matQ_{\le i} \, \left\{ \matS_{i}  - \matQ_{\le i}^\top [\Delta A]^{\langle i \rangle}   \matQ_{\ge i+1} \right\}\, \matQ_{\ge i+1}^\top.
\]     
These results are furthermore valid for any ordering of the initial value problems.
\end{theorem}
\begin{proof} First, observe that each $P_i^+$ and $P_i^-$ maps onto the current tangent space of $\M$ and that $Y_0 \in \M$. Hence, each $Y_i^+(t)$ and $Y_i^-(t)$ will stay on $\M$. % for $t_1-t_0$ sufficiently small. 
We may therefore assume that $Y_i^+(t)$ and $Y_i^-(t)$ admit TT/MPS decompositions of equal TT/MPS rank for $t \in [t_0,t_1]$.

By writing $Y_i^+(t)$ in a time-dependent recursive SVD,
\[
 [Y_i^+(t)]^{\langle i \rangle} = 
 (\matI \otimes \matQ_{\le i-1}(t)) \matQ_{i}^<(t) \matS_{i}(t) \matQ_{\ge i+1}^\top(t),
\]
we see from \eqref{eq:def_Pi_plus} that
\[
 [P_i^+(\dt A)] ^{\langle i \rangle} =  
 (\matI \otimes \matQ_{\le i-1}(t) \matQ^\top_{\le i-1}(t)) \,[ \dt A \,]^{\langle i \rangle} \, \matQ_{\ge i+1}(t) \matQ_{\ge i+1}^\top(t).
\]
Hence the differential equation  $\dt Y_i^+ = P_i^+(\dt A)$ implies 
\begin{align*}
  &(\matI \otimes \dt\matQ_{\le i-1}(t)) \matQ_{i}^<(t) \matS_{i}(t) \matQ_{\ge i+1}^\top(t) +
  (\matI \otimes \matQ_{\le i-1}(t)) \frac{\ud}{\ud t} [\matQ_{i}^<(t) \matS_{i}(t)]\ \matQ_{\ge i+1}^\top(t) \\
  &\qquad+  (\matI \otimes \matQ_{\le i-1}(t)) \matQ_{i}^<(t) \matS_{i}(t) \dt\matQ_{\ge i+1}^\top(t) \\
  &= (\matI \otimes \matQ_{\le i-1}(t)) (\matI \otimes \matQ^\top_{\le i-1}(t)) [ \dt A \,]^{\langle i \rangle}  \matQ_{\ge i+1}(t) \matQ_{\ge i+1}^\top(t).
\end{align*}
By choosing $\dt\matQ_{\le i-1}(t) = 0$ and $\dt\matQ_{\ge i+1}(t) = 0$, the above identity is satisfied when
\begin{align}\label{eq:dK_dt}
    &\frac{\ud}{\ud t} [\matQ_{i}^<(t) \matS_{i}(t)] = (\matI \otimes \matQ^\top_{\le i-1}(t)) [ \dt A \,]^{\langle i \rangle}  \matQ_{\ge i+1}(t).
\end{align}
Using the initial condition $Y_i^+(0)$, the solution of these differential equations becomes
\begin{align*}
    &\matQ_{\le i-1}(t) = \matQ_{\le i-1}(0), \qquad \matQ_{\ge i+1}(t) = \matQ_{\ge i+1}(0), \\    
    &\matQ_{i}^<(t) \matS_{i}(t) = \matQ_{i}^<(0) \matS_{i}(0) +  (\matI \otimes \matQ^\top_{\le i-1}(0)) [ A(t) - A(0)]^{\langle i \rangle}  \matQ_{\ge i+1}(0),
\end{align*}
which proves the statement for $[Y_i^+(t_1)]^{\langle i \rangle}$. Now, writing 
\begin{align*}
 [Y_i^+(t_1)]^{\langle i \rangle}  &= (\matI \otimes \matQ_{\le i-1}) \, \matQ_{i}^< \matS_{i}  \, \matQ_{\ge i+1}^\top + 
 (\matI \otimes \matQ_{\le i-1} \matQ_{\le i-1}^\top) \,  [\Delta A]^{\langle i \rangle}   \matQ_{\ge i+1}  \, \matQ_{\ge i+1}^\top  \\
 &= [Y_i^+(0)]^{\langle i \rangle} + (\matI \otimes \matP_{\le i-1}) [\Delta A]^{\langle i \rangle}  \matP_{\ge i+1},
\end{align*}
we have also proven the first statement of the theorem.

Since the previous derivation is valid for any initial condition, it does not depend on a specific ordering of the initial value problems. The derivation for $Y_i^-(t_1)$ is analogous to that of $Y_i^+(t_1)$.
\end{proof}     

In a similar way as for the proof of Corollary~\ref{cor:eq_PX}, one can show that the projector~\eqref{eq:def_Pi_plus} also satisfies
\begin{align}
 P_i^+(Z) &= \ten_{i-1} \bigl[ \matP_{\le i-1} \matZ^{\langle i-1 \rangle} (\matP_{\ge i+1} \otimes \matI_{n_i})  \bigr],  &  (1 \leq i \leq d). \label{eq:def_Pi_plus_b}
\end{align}
This definition is useful when $Y_i^+(t_0)$ is given as (see \S\ref{sec:full_sweep})
\begin{equation*}
 [Y_i^+(t_0)]^{\langle i-1 \rangle}  = \matQ_{\le i-1} \matS_{i-1} \matQ_{i}^{>\top} (\matQ_{\ge i+1}^\top \otimes \matI_{n_i}).
\end{equation*}
In that case, we have
\begin{equation}\label{eq:Y_i_b}
 [Y_i^+(t_1)]^{\langle i-1 \rangle}  = \matQ_{\le i-1} \, \left\{ \matS_{i-1} \matQ_{i}^{>\top} \ + \matQ_{\le i-1}^\top [\Delta A]^{\langle i-1   \rangle}   (\matQ_{\ge i+1} \otimes \matI_{n_i}) \right\} \, (\matQ_{\ge i+1}^\top \otimes \matI_{n_i}).
\end{equation}

\subsection{Efficient implementation as a sweeping algorithm}\label{sec:sweep_forward}

\begin{figure}
     \centering
          \vspace{0.2cm}
     \def\svgwidth{\columnwidth}
     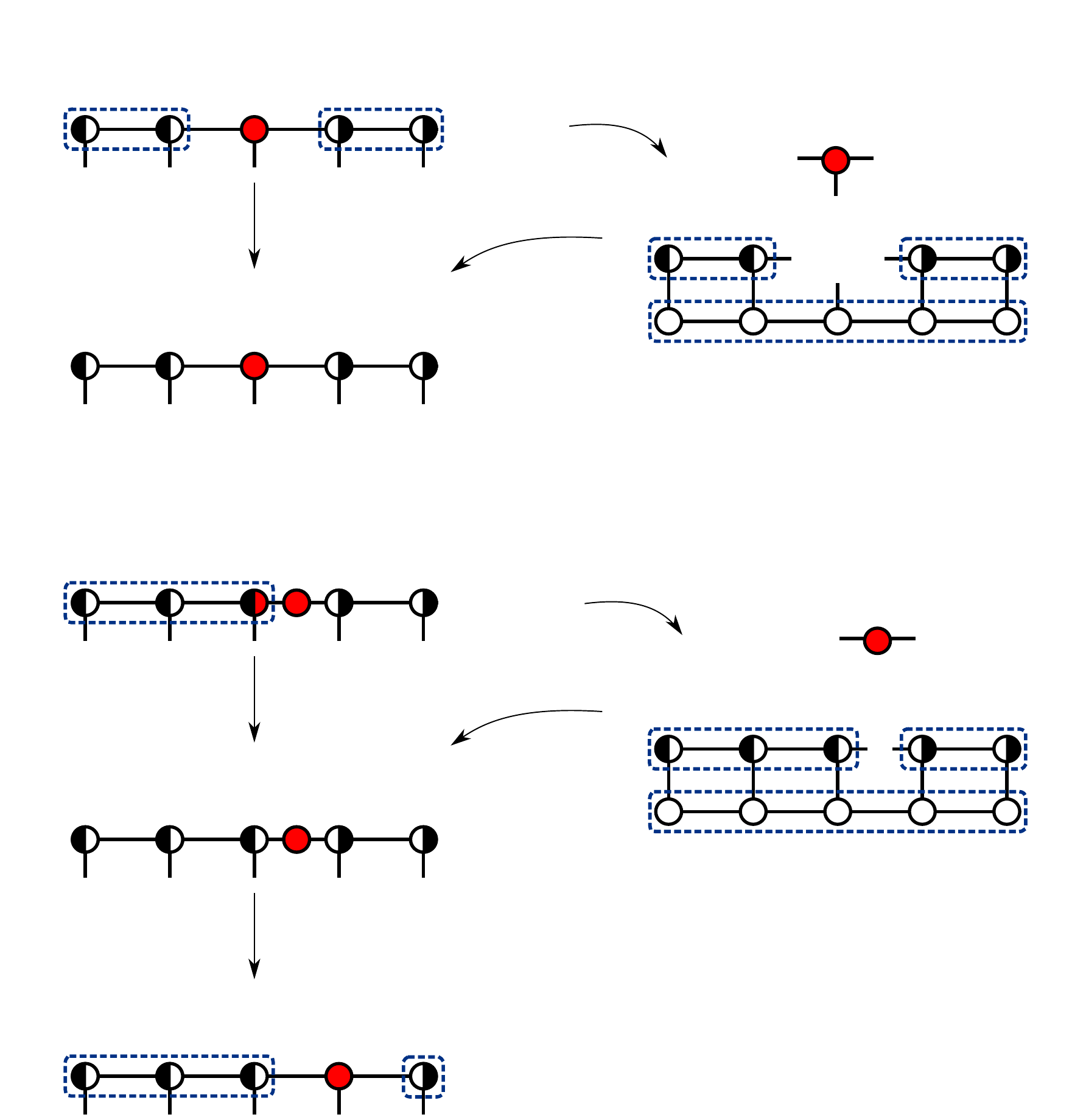
     \vspace{0.2cm}
     \caption{The two sweeping algorithms update the cores selectively throughout the time stepping computations. Shown for the forward sweep when computing $Y_i^+$ and $Y_i^-$.}\label{fig:sweep}
 \end{figure}

Theorem~\ref{thm:sol_split} can be turned into an efficient scheme by updating the cores of the tensor $Y(t_0)$ from left to right. Our explanation will be high level, focusing only on pointing out which cores stay constant and which need to be updated throughout the sweep. A graphical depiction of the resulting procedure using tensor networks is given in Fig.~\ref{fig:sweep}. More detailed implementation issues are deferred to \S\ref{sec:impl_details}.

{\it Preparation of $Y_0$.\/} Before solving the substeps, we prepare the starting value $Y_0$ as follows. Write $Y=Y_0$ for notational convenience and suppose
$$
\matY^{\langle 1 \rangle}(t_0) = \matY_{\le 1}(t_0) \, \matY_{\ge 2}^\top(t_0).
$$ 
By orthogonalization from the right we decompose
$ \matY_{\ge 2}(t_0)=  \matQ_{\ge 2}(t_0)  \matR_{2}(t_0)$, 
%and from the left we factorize $\matY_{\le 1}(t_0)=\matC_1^<(t_0) = \matQ_{\le 1}(t_0)  \matR_{\le 1}(t_0)$,
so that we obtain the right-orthogonalized factorization
$$
\matY^{\langle 1 \rangle}(t_0) = \matK^<_1(t_0) \matQ_{\ge 2}^\top(t_0) 
$$
with $\matK_1^<(t_0)= \matY_{\le 1}(t_0) \matR_{2}^\top(t_0) \in \R^{n_1\times r_1}$ the first core of $Y(t_0)$.

{\it Computation of $Y_1^+$.\/} Denote $Y=Y_1^+$. Since $ \matQ_{\ge 2}^\top(t_0)  \matQ_{\ge 2}(t_0) = \matI_{r_1}$, we have  that $\matP_{\ge 2}(t_0)=
\matQ_{\ge 2}(t_0)  \matQ_{\ge 2}^\top (t_0)$. Applying Theorem~\ref{thm:sol_split} gives
\[
\matY^{\langle 1 \rangle}(t_1) = \matK_1^<(t_1) \, \matQ_{\ge 2}^\top(t_0), 
\]
with 
\[
\matK_1^<(t_1) = \matK_{1}^<(t_0) + 
 \bigl(\matA^{\!\langle 1 \rangle}(t_1)-\matA^{\!\langle 1 \rangle}(t_0)\bigr) \matQ_{\ge 2}(t_0).
 \]
Observe that compared to $Y(t_0)$ only the first core $K_1(t_1)$ of $Y(t_1)$ is changed, while all the others (that is, those that make up $\matQ_{\ge 2}(t_0)$) stay constant.
Hence, after computing the QR decomposition
 $$
 \matK_1^<(t_1) = \matQ_{1}^<(t_1)  \matR_1(t_1),
 $$ 
we obtain a recursive SVD for $Y(t_1) = Y_1^+(t_1)$,
 $$
\matY^{\langle 1 \rangle}(t_1) = \matQ_{\le 1}(t_1) \matR_1(t_1) \matQ_{\ge 2}^\top(t_0)  \qquad \text{with} \quad \matQ_{\le 1}(t_1)  = \matQ_{1}^<(t_1) .
$$

{\it Computation of $Y_i^-$ with $i=1,\ldots,d-1$.\/} The computation for $Y_1^-$ follows the same pattern as for arbitrary $Y_i^-$, so we explain it directly for $Y_i^-$.

We require that the initial value $Y_i^-(t_0) = Y_{i}^+(t_1)$ is available as a recursive SVD in node $i$. This is obviously true for $Y_1^+(t_1)$ and one can verify by induction that it is also true for $Y_i^+(t_1)$ with $i>1$, whose computation is explained below. Denoting $Y = Y_i^-$, we have in particular
 $$
\matY^{\langle i \rangle}(t_0) = \matQ_{\le i}(t_1) \matR_i(t_1) \matQ_{\ge i+1}^\top(t_0),
$$
with $\matQ_{\le i}^\top(t_1) \matQ_{\le i}(t_1) = \matI_{r_i}=\matQ_{\ge i+1}^\top(t_0) \matQ_{\ge i+1}(t_0)$. This means we can directly apply Theorem~\ref{thm:sol_split} for the computation of $Y(t_1)$ and obtain
\begin{equation}\label{eq:Yi_min}
\matY^{\langle i     \rangle} (t_1) = \matQ_{\le i}(t_1) \matS_i(t_1) \matQ_{\ge i+1}^\top(t_0),
\end{equation}
where $\matS_i(t_1) \in \R^{r_i\times r_i}$ is given as
\begin{equation}\label{eq:Yi_min_Si}
\matS_i(t_1) = \matR_i(t_1) -  \matQ_{\le i}^\top(t_1)
\bigl( \matA^{\!\langle i \rangle}(t)-\matA^{\!\langle i \rangle}(t_0) \bigr)
 \matQ_{\ge i+1}(t_0).
\end{equation}
Observe that we maintain a recursive SVD in $i$ for  $Y_i^-(t_1)$ without having to orthogonalize the matrices $\matQ_{\le i}(t_1)$ or $\matQ_{\ge i+1}(t_0)$.

 {\it Computation of $Y_i^+$ with $i=2,\ldots,d$.\/}   In this case, the initial value $Y_i^+(t_0) = Y_{i-1}^-(t_1)$ is available as a recursive SVD in node $i-1$. Denoting $Y = Y_i^+$, then it is easily verified by induction that
$$
\matY^{\langle i-1\rangle}(t_0) = \matQ_{\le i-1}(t_1)\, \matS_{i-1}(t_1) \, \matQ_{\ge i}^\top(t_0),
$$
with $\matQ_{\le i-1}^\top(t_1) \matQ_{\le i-1}(t_1) = \matI_{r_{i-1}} = \matQ_{\ge i}^\top(t_0) \matQ_{\ge i}(t_0)$. Recalling the relations \eqref{eq:rec_svd_i} and~\eqref{eq:rec_svd_iplus}, we can transform this $(i-1)$th unfolding into the $i$th unfolding, 
\begin{equation}\label{eq:Yi_plus_t0}
\matY^{\langle i \rangle}(t_0) = (\matI_{n_i}\otimes \matQ_{\le i-1}(t_1))\, \matK_i^<(t_0) \,  \matQ_{\ge i+1}^\top(t_0)
\end{equation}
where $\matK_i^<(t_0)=(\matI_{n_i}\otimes  \matS_{i-1}(t_1)) \matQ_{i}^<(t_0)$ is the left unfolding of the $i$th core $K_i(t_0)$ of $Y(t_0)$. Theorem~\ref{thm:sol_split} then leads to
\begin{equation}\label{eq:Yi_plus_t1}
\matY^{\langle i \rangle}(t_1) = (\matI_{n_i}\otimes \matQ_{\le i-1}(t_1))\matK_i^<(t_1) \, \matQ_{\ge i+1}^\top(t_0),
\end{equation}
where $\matK_i^<(t_1) \in \R^{(r_{i-1}n_i)\times r_i}$ is given by
$$
 \matK_i^<(t_1) = \matK_i^<(t_0) + (\matI_{n_i}\otimes \matQ_{\le i-1}^\top(t_1))
 \bigl(\matA^{\!\langle i \rangle}(t_1)-\matA^{\!\langle i \rangle}(t_0)\bigr) \matQ_{\ge i+1}(t_0).
 $$
 Since now only the $i$th core $K_i(t_1)$ of $Y(t_1)$ has changed, one QR decomposition
 \begin{equation}\label{eq:Yi_plus_QR}
 \matK_i^<(t_1) = \matQ_{i}^<(t_1)  \matR_i(t_1),
 \end{equation}
 suffices to obtain a recursive SVD of $Y_i^+(t_1) = Y(t_1)$ at node $i$,
 $$
\matY^{\langle i \rangle}(t_1) =  \matQ_{\le i}(t_1) \matR_i(t_1) \matQ_{\ge i+1}^\top(t_0), \qquad \text{with} \quad \matQ_{\le i}(t_1) = (\matI_{n_i}\otimes \matQ_{\le i-1}(t_1)) \matQ_{i}^<(t_1). 
$$

 {\it Next time step.\/} The final step $Y_d^+(t_1)$ will be an approximation to $Y(t_1)$ and consists of a left-orthogonal $\matQ_{\le d}(t_1)$. If we now want to continue with the time stepper to approximate $Y(t_2)$ for $t_2 > t_1$, we need to apply the scheme again using $Y_d^+(t_1)$ as initial value. This requires a new orthogonalization procedure from right to left, since the initial value for the sweep has to be right orthogonalized.

\subsection{Second-order scheme by a back-and-forth sweep}\label{sec:full_sweep}

% The final step of the split projector will have computed $Y_d^+$ as
%   $$
%   \matY^{\langle d \rangle}(t_1) = (\matI_{n_d}\otimes \matQ_{\le d-1}(t_1))\matK_d^<(t_1), 
%   $$
% where
%   $$
%    \matK_d^<(t_1) = \matK_d^<(t_0) + (\matI_{n_d}\otimes \matQ_{\le d-1}^\top(t_1))
%    \bigl(\matA^{\!\langle d \rangle}(t)-\matA^{\!\langle d \rangle}(t_0)\bigr) 
%    $$
% with
%   \[
%   \matK_d^<(t_1) =(\matI_{n_d}\otimes  \matS_{d-1}(t_1)) \matQ_{d}^<(t_0) \in \R^{(r_{d-1}n_d)\times 1}.
%   \]
% In this way all cores will have been updated.

In many cases, it is advisable to compose the scheme from above with its adjoint instead of only orthogonalizing and continuing with the next step. In particular, the Strang splitting consists of first computing the original splitting scheme on $t \in [t_0, \th]$ with $\th = (t_1+t_0)/2$ and then applying the adjoint of this scheme on $t \in [\th, t_1]$.  The result will be a symmetric time stepper of order two; see, e.g., \cite[II.5]{Hairer:2006}.

For our splitting, the adjoint step is simply solving the split differential equations in reverse order. Since Theorem~\ref{thm:sol_split} is independent of the ordering of the differential equations, we can again use its closed-form solutions to derive an efficient sweeping algorithm for this adjoint step. We briefly explain the first three steps and refer to  Algorithm~\ref{al:int} for the full second-order scheme. Observe that this scheme can be seen as a full back-and-forth sweep.

Denote the final step of the forward sweep on $t \in [t_0, \th]$ by $\widehat Y = Y_d^+(\th)$. It satisfies (recall that $t_1$ takes the role of $\th$ in the derivations above)
\[
\widehat\matY^{\langle d \rangle}(\th) = (\matI_{n_d}\otimes \matQ_{\le d-1}(\th))\matK_d^<(\th).
\]
with
\[
  \matK_d^<(\th) = \matK_d^<(t_0) + (\matI_{n_d}\otimes \matQ_{\le d-1}^\top(\th))
  \bigl(\matA^{\!\langle d \rangle}(\th)-\matA^{\!\langle d \rangle}(t_0)\bigr).
\]
The first substep of the adjoint scheme consists of solving 
\begin{align*}
 \dt Y_d^+ &= P_d^+( \dt A), \qquad Y_d^+(\th) = \widehat Y,
\end{align*} 
for $t \in [\th,t_1]$. Denote $Y=Y_d^+$. We can directly apply Theorem~\ref{thm:sol_split} to obtain
\[
\matY^{\langle d \rangle}(t_1) = (\matI_{n_d}\otimes \matQ_{\le d-1}(\th))\matK_d^<(t_1)
\]
with
\begin{align*}
  \matK_d^<(t_1) &= \matK_d^<(\th) + (\matI_{n_d}\otimes \matQ_{\le d-1}^\top(\th))
  \bigl(\matA^{\!\langle d \rangle}(t_1)-\matA^{\!\langle d \rangle}(\th)\bigr) \\
  &= \matK_d^<(t_0) + (\matI_{n_d}\otimes \matQ_{\le d-1}^\top(\th))
    \bigl(\matA^{\!\langle d \rangle}(t_1)-\matA^{\!\langle d \rangle}(t_0)\bigr) .
\end{align*}
Hence, the last substep of the forward sweep and the first of the backward sweep can be combined into one.

The second substep amounts to solving
\begin{align*}
 \dt Y_{d-1}^- &= -P_{d-1}^-( \dt A), \qquad Y_{d-1}^-(\th) = Y_{d}^+(t_1).
\end{align*} 
Let $Y=Y_{d-1}^-$. Then we can write the initial condition as
\[
 \matY^{\langle d-1 \rangle}(\th) =  \matQ_{\le d-1}(\th) \matK_d^>(t_1)^\top
\]
and using the QR decomposition $\matK_d^>(t_1) = \matQ_{d}^>(t_1) \matR_{d-1}(t_1)$ also as
\[
 \matY^{\langle d-1 \rangle}(\th) =  \matQ_{\le d-1}(\th) \matR^\top_{d-1}(t_1) \matQ_{\ge d}^{\top}(t_1),
\]
where $\matQ_{\ge d}(t_1) = \matQ_{d}^>(t_1)$. Applying Theorem~\ref{thm:sol_split}, we obtain
$$
\matY^{\langle d-1     \rangle} (t_1) = \matQ_{\le d-1}(\th) \matS_{d-1}^\top(t_1) \matQ_{\ge d}^\top(t_1),
$$
where
\[
\matS_{d-1}^\top(t_1) = \matR_{d-1}^\top(t_1) -  \matQ_{\le d-1}^\top(\th)
\bigl( \matA^{\!\langle d-1 \rangle}(t_1)-\matA^{\!\langle d-1 \rangle}(\th) \bigr)
 \matQ_{\ge d}(t_1).
\]

For the third substep 
\begin{align*}
 \dt Y_{d-1}^+ &= P_{d-1}^-( \dt A), \qquad Y_{d-1}^+(\th) = Y_{d-1}^-(t_1),
\end{align*} 
we denote $Y = Y_{d-1}^+$. In this case, unfold using \eqref{eq:rec_svd_i_b} and \eqref{eq:rec_svd_imin} the computed quantity $Y_{d-1}(t_1)$ from above  as
\[
 \matY^{\langle d-2 \rangle}(t_1) =  \matQ_{\le d-2}(\th) \matK^{>\top}_{d-1}(\th) (\matQ_{\ge d}^{\top}(t_1) \otimes \matI_{n_{d-1}}),
\]
with $\matK_{d-1}^{>\top}  (\th) = \matQ_{d-1}^{>\top}(\th) (\matS^\top_{d-1}(t_1) \otimes \matI_{n_{d-1}})$. From here on, all subsequent computations are straightforward if we use \eqref{eq:Y_i_b} to compute $Y_i^+(t_1)$.

\begin{algorithm} \label{al:int}
\DontPrintSemicolon
\SetKwComment{tcp}{$\triangleright$~}{}
\KwData{$K_1(t_0)\in\R^{r_0\times n_1 \times r_1}$, $Q_i(t_0)\in\R^{r_{i-1}\times n_i \times r_i}$ with $\matQ_i^{>\top}(t_0)\matQ_i^>(t_0)=\matI_{r_{i-1}}$ for $i=2,\dots,d$; $t_0, t_1$ }
\KwResult{$K_1(t_1)\in\R^{r_0\times n_1 \times r_1}$, $Q_i(t_1)\in\R^{r_{i-1}\times n_i \times r_i}$ with $\matQ_i^{>\top}(t_1)\matQ_i^>(t_1)=\matI_{r_{i-1}}$ for $i=2,\dots,d$}
\Begin{
\nl set $\th =(t_0+t_1)/2$. \tcp*[r]{Initialization}

\nl set $\Dl = A(\th) - A(t_0)$ and $\Dr  = A(t_1) - A(\th)$.

% comment: start forward sweep

\nl set $\matK_1^<(\th) = \matK_1^<(t_0) +   
\matDl^{\!\langle 1 \rangle} \matQ_{\ge 2}(t_0)$.  \tcp*[r]{Forward}
 
\nl compute QR factorization $\matK_1^<(\th) = \matQ_1^<(\th)\matR_1(\th)$.

\nl set $\matQ_{\le 1}(\th) = \matQ_1^<(\th)$.

\nl set $ \matS_1(\th) = \matR_1(\th) - \matQ_{\le 1}^\top(\th)
 \matDl^{\!\langle 1 \rangle} \matQ_{\ge 2}(t_0)$.

\nl \For{$i=2$ \KwTo $d-1$}
{ 
\nl
{set $\matK_i^<(t_0)=(\matI_{n_i}\otimes \matS_{i-1}(\th))\matQ_i^<(t_0)$}.

\nl set $\matK_i^<(\th) = \matK_i^<(t_0) + (\matI_{n_i} \otimes \matQ^\top_{\le i-1}(\th))
 \matDl^{\!\langle i \rangle} \matQ_{\ge i+1}(t_0)$.
 
\nl compute QR factorization $\matK_i^<(\th) = \matQ_i^<(\th)\matR_i(\th)$

\nl set $\matQ_{\le i}(\th) = (\matI_{n_i} \otimes \matQ_{\le i-1}(\th)) \matQ_i^<(\th)$

\nl set $ \matS_i(\th) = \matR_i(\th) - \matQ^\top_{\le i}(\th)
 \matDl^{\!\langle i \rangle} \matQ_{\ge i+1}(t_0)$
}

\nl set $\matK_d^<(t_1) = \matK_d^<(t_0) +
 (\matI_{n_d} \otimes \matQ^\top_{\le d-1}(\th))
 \bigl(\matDl^{\!\langle d \rangle} + \matDr^{\!\langle d \rangle}\bigr)$   \tcp*[r]{Backward}

 \nl compute QR factorization $\matK_d^>(t_1) = \matQ_d^>(t_1)\matR_{d-1}(t_1)$
 
 \nl set $\matQ_{\ge d}(t_1) = \matQ_d^>(t_1)$
 
 \nl set $\matS_{d-1}(t_1) = \matR_{d-1}(t_1) - \matQ_{\ge d}^\top(t_1)  (\matDr^{\!\langle d-1 \rangle})^\top \matQ_{\le d-1}(\th)$
 
 \nl \For{$i=d-1$ {\bf down to} $2$}{
% \nl {\bf if $i<d$ then} 
\nl set
 $\matK_i^>(\th) = (\matS_{i}(t_1) \otimes \matI_{n_i} ) \matQ_{i}^>(\th)$

 \nl set $\matK_i^>(t_1) =  \matK_i^>(\th) + (\matQ^\top_{\ge i+1}(t_1) \otimes  \matI_{n_i})
  (\matDr^{\!\langle i-1 \rangle})^\top \matQ_{\le i-1}(\th)$
  
  \nl compute QR factorization $\matK_i^>(t_1) = \matQ_i^>(t_1)\matR_{i-1}(t_1)$
  
  \nl set $\matQ_{\ge i}(t_1) = (\matQ_{\ge i+1}(t_1) \otimes   \matI_{n_i}) \matQ_i^>(t_1)$
 
 \nl set $\matS_{i-1}(t_1) = \matR_{i-1}(t_1) - \matQ_{\ge i}^\top(t_1)  (\matDr^{\!\langle i-1 \rangle})^\top  \matQ_{\le i-1}(\th)$
 }
 
\nl  set $\matK_1^>(\th) = (\matS_{1}(t_1) \otimes   \matI_{n_1}) \matQ_{1}^>(\th)$
 
  \nl set $\matK_1^>(t_1) =  \matK_1^>(\th) + (\matQ^\top_{\ge 2}(t_1) \otimes  \matI_{n_1}) (\matDr^{\!\langle 0 \rangle})^\top$

}
\caption{Step of the split projector integrator of second order}
\end{algorithm}
\smallskip

\section{Exactness property of the integrator}

We show that the splitting integrator is exact when $A(t)$ is a tensor of constant TT/MPS rank $\TTrank$. This is similar to Theorem 4.1 in \cite{Lubich:2014} for the matrix case, except that in our case we require the rank of  $A(t)$ to be exactly $\TTrank$ and not merely bounded by $\TTrank$. Note, however, that the positive singular values of unfoldings of $A(t)$ can be arbitrarily small. \endCL

\begin{theorem}\label{thm:exactness}    
    Suppose $A(t) \in \M$ for $t \in [t_0,t_1]$. Then, for sufficiently small $t_1-t_0>0$ the splitting  integrators of orders one and two are exact when started from $Y_0  = A(t_0)$. For example, $Y_d^+(t_1) = A(t_1)$ for the first-order integrator.
\end{theorem}    

The proof of this theorem follows trivially from the following lemma.

\begin{lemma}\label{lem:steps_exactness}    
Suppose $A(t) \in \M$ for $t \in [t_0,t_1]$ with recursive SVDs
\[
 [\matA(t)]^{\langle i \rangle} = \matQ_{\le i}(t) \, \matS_{i}(t) \, \matQ_{\ge i+1}^\top(t) \qquad \text{for $i=0,1,\ldots,d$}.
\]
Let $Y_0 = A(t_0)$, then for sufficiently small $t_1-t_0>0$ the consecutive steps in the splitting integrator of \S\ref{sec:forward_sweep} satisfy
\[
 Y_i^+(t_1) = P_{\ge i+1}^{(0)} A(t_1) \qquad \text{and} \qquad  Y_i^-(t_1) = P_{\le i}^{(1)} A(t_0) \qquad \text{for $i=1,2,\ldots,d$},
\]
where 
\begin{align*}
 &P_{\ge i+1}^{(0)} Z =  \ten_{i} \bigl(   \matZ^{\langle i \rangle} \matQ_{\ge i+1}(t_0)\matQ^\top_{\ge i+1}(t_0) \bigr),\\
 &P_{\le i}^{(1)} Z = \ten_{i} \bigl( \matQ_{\le i}(t_1)\matQ^\top_{\le i}(t_1) \matZ^{\langle i \rangle} \bigr)  .
\end{align*}
\end{lemma}

Before proving this lemma, we point out that the assumption of sufficiently small $t_1-t_0$ is only because the matrices $\matQ^\top_{\le i-1}(t_1) \matQ_{\le i-1}(t_0)$ and $\matQ_{\ge i}^\top(t_1) \matQ_{\ge i}(t_0)$ need to be invertible. Since the full column-rank matrices $\matQ_{\le i}(t)$ and $\matQ_{\ge i}(t)$ can be chosen continuous functions in $t$, this is always satisfied for $t_1-t_0$ sufficiently small. It may however also hold for larger values of $t_1-t_0$.

% In addition, it also holds almost surely for any $\delta$ if we take $A(t)$ a generic curve in $\M$.\footnote{While the smallest singular value of $\matQ^\top_{\le i-1}(\delta) \matQ_{\le i-1}(0)$ and $\matQ_{\ge i}^\top(\delta) \matQ_{\ge i}(0)$ can concentrate around zero, the probability that the matrices are rank deficient are zero a.s. Since Lemma~\ref{lem:steps_exactness} is not about numerical since}

\begin{proof}
    The proof proceeds by induction on $i$ from left to right. Since $Y_1^+(t_0) = A(t_0)$, we can include the case for $Y_1^+(t_1)$ in our proof below for general $i$ by putting $Y_{0}^-(t_1) = Y_0^+(t_0)$ and $P_{\le 0}^{(1)} = 1$. 
    
    Now, suppose the statement to be true for $i > 1$. Then, $Y_i^+(t_0) = Y_{i-1}^-(t_1) = P_{\le i-1}^{(1)} A(t_0)$, which gives
    \begin{align*}
     [\matY_i^+(t_0)]^{\langle i-1 \rangle} &= \matQ_{\le i-1}(t_1)\, \matQ^\top_{\le i-1}(t_1) \matQ_{\le i-1}(t_0)  \matS_{i-1}(t_0) \, \matQ_{\ge i}^\top(t_0) \\
      &= \matQ_{\le i-1}(t_1)\, \matS_{i-1}^+ \, \matQ_{\ge i}^\top(t_0).
    \end{align*}
    Observe that $Y_i^+(t_0)\in \M$ since  $\matS_{i-1}^+ = \matQ^\top_{\le i-1}(t_1) \matQ_{\le i-1}(t_0)  \matS_{i-1}(t_0)$ is of full rank for $t_1 - t_0$ sufficiently small. Hence, from \eqref{eq:rec_svd_i}--\eqref{eq:rec_svd_iplus} we obtain
    \[
     [\matY_i^+(t_0)]^{\langle i \rangle} = (\matI_{n_i} \otimes \matQ_{\le i-1}(t_1))\ (\matI_{n_i} \otimes \matS_{i-1}^+) \ \matQ_{i}^{<\top}(t_0) \ \matQ_{\ge i+1}^\top(t_0).
    \]
Comparing to \eqref{eq:def_Pi_plus}, we see that the projector onto the tangent space at $Y_i^+(t_0)$ equals $P_i^+ = P_{\le i-1}^{(1)} \, P_{\ge i+1}^{(0)} = P_{\ge i+1}^{(0)} P_{\le i-1}^{(1)} $. The  previous identities give with Theorem~\ref{thm:sol_split} that
\begin{align*}
 Y_i^+(t_1) &= Y_i^+(t_0) + P_i^+ A(t_1) - P_i^+ A(t_0) \\
 &= P_{\le i-1}^{(1)} A(t_0) + P_{\ge i+1}^{(0)} \,  P_{\le i-1}^{(1)} A(t_1) - P_{\le i-1}^{(1)} \, P_{\ge i+1}^{(0)} A(t_0) = P_{\ge i+1}^{(0)} A(t_1), 
\end{align*}
where we used $P_{\ge i+1}^{(0)} A(t_0) = A(t_0)$ and $P_{\le i-1}^{(1)} A(t_1) = A(t_1)$. 

Continuing with $Y_i^-(t_0) = Y_{i}^+(t_1) = P_{\ge i+1}^{(0)} A(t_1)$, we have
\begin{align*}
 [Y_i^-(t_0)]^{\langle i \rangle} &= \matQ_{\le i}(t_1)\, \matS_{i}(t_1) \, \matQ_{\ge i+1}^\top(t_1) \matQ_{\ge i+1}(t_0) \matQ_{\ge i+1}^\top(t_0) \\
  &= \matQ_{\le i}(t_1)\, \matS_{i}^- \, \matQ_{\ge i+1}^\top(t_0). 
\end{align*}
This is again a recursive SVD with full rank $\matS_i^- = \matS_{i}(t_1) \, \matQ_{\ge i+1}^\top(t_1) \matQ_{\ge i+1}(t_0)$. Comparing to \eqref{eq:def_Pi_min}, we have $P_i^- = P_{\le i}^{(1)} P_{\ge i+1}^{(0)} = P_{\ge i+1}^{(0)} P_{\le i}^{(1)} $ and Theorem~\ref{thm:sol_split} gives
\begin{align*}
 Y_i^-(t_1) &= Y_i^-(t_0) - P_i^- A(t_1) + P_i^- A(t_0) \\
 &= P_{\ge i+1}^{(0)} A(t_1) -  P_{\ge i+1}^{(0)} P_{\le i}^{(1)} A(t_1) + P_{\le i}^{(1)} P_{\ge i+1}^{(0)} A(t_0) = P_{\le i}^{(1)}  A(t_0),
\end{align*}
where we used $P_{\le i}^{(1)} A(t_1) = A(t_1)$.  This concludes the proof.
\end{proof}    

Now, Theorem~\ref{thm:exactness} is a simple corollary.

{\em Proof of Theorem~\ref{thm:exactness}.}
For the forward sweep (that is, the first-order scheme), Lemma~\ref{lem:steps_exactness} immediately gives exactness since $Y_d^+(t_1) = P_{\ge d+1}^{(0)} A(t_1) = A(t_1)$ with $\matQ_{\ge d+1}(t_0)=1$. The second-order scheme composes this forward sweep with a backward sweep involving the same substeps. It is not difficult to prove the analogous version of Lemma~\ref{lem:steps_exactness} for such a backward ordering such that we establish exactness for the second-order scheme too. \endproof

\section{Numerical implementation and experiments}

We consider two numerical experiments. First, we use the splitting integrator for the integration of a time-dependent molecular Schr\"odinger
equation with a model potential. %This is one of the most promising
%applications for the new integrator, since the whole idea of dynamical
%low-rank approximation was proposed in this community. 
In the second experiment, we use one step of the splitting integrator as a retraction on the manifold of TT/MPS tensors and perform a
Newton--Schultz iteration for approximate matrix inversion.

\subsection{Implementation details}\label{sec:impl_details}

% Starting form a tensor train orthogonalized from the right, a step of the projector-splitting integrator described in the previous section yields an updated tensor train that is orthogonalized from the left. Instead of restarting the algorithm after a reorthogonalization from the right, it is preferable to add another step with the reversed order of the splitting, starting from the obtained tensor train orthogonalized from the left and yielding a tensor train that is again orthogonalized from the right. The combined method with a forward and a backward sweep is a symmetric, Strang-type splitting and yields an integrator of order 2.

As explained in \S\ref{sec:sweep_forward}--\ref{sec:full_sweep}, the integrator updates the cores $K_i$ and matrices $\matS_i$ in a forward, and possibly, backward ordering. Except for the (relatively cheap) orthogonalizations of the cores, the most computationally intensive part of the algorithm is computing these updates. For example, in the forward sweep, we need to compute the contractions (see Fig.~\ref{fig:sweep})
\begin{align*}
 \Delta_i^+ &=  (\matI \otimes \matQ^\top_{\le i-1}(t_1)) \, [ A(t_1) - A(t_0) \,]^{\langle i \rangle}  \,  \matQ_{\ge i+1}(t_0),\\
  \Delta_i^- &= \matQ^\top_{\le i}(t_1) \,[ A(t_1) - A(t_0) \,]^{\langle i \rangle}  \,  \matQ_{\ge i+1}(t_0).
\end{align*}
It is highly recommended to avoid constructing the matrices $\matQ_{\le i}$ and $\matQ_{\ge i}$ explicitly when computing $\Delta_i^+,\Delta_i^-$ and instead exploit their TT/MPS structure. How this can be done, depends mostly on the structure of the increments $A(t_1) - A(t_0)$. In particular, the contractions are computed inexpensively if $A(t)$ is itself a linear combination of  TT/MPS tensors, possibly having different rank than $Y$, and a sparse tensor.

%This is for example the case for Laplacian operators or more general Hamiltonians with short-range interactions; see, e.g., \cite{Kressner:2013}.

% The update formulas for $K_i$ and $\matS_i$ involve contractions of the increments
% $A(\th)-A(t_0)$ and $A(t_1)-A(\th)$ with the tensors $Q_{\le i}$ and $Q_{\ge i+1}$ that are constructed from the 3-tensors $Q_i$ recursively as 
% $
% \matQ_{\le i}=(\matI_{n_i}\otimes \matQ_{\le i-1}) \matQ_i^<$ and
% $\matQ_{\ge i+1}=(\matI_{n_i}\otimes \matQ_{\ge i+2}) \matQ_{i+1}^>$.

The computation of $K_i$ and $\matS_i$ changes when the tensor $A(t)$ is not given explicitly, but  determined as the solution of a tensor differential equation 
\[
 \dt A(t) = F(t,A(t)).
\]
In case of a forward sweep, $Y_i^+(t_1)$ is obtained as the evaluation at $t=t_1$ of
\[
 Y_i^+(t) = (\matI \otimes \matQ_{\le i-1}(t_1)) \, \matK_i^<(t) \,  \matQ_{\ge i+1}^\top(t_0),
\]
where  $\matK_i^<(t)=\matQ_{i}^<(t) \matS_{i}(t)$ satisfies \eqref{eq:dK_dt}. Hence, for $\dt A(t) = F(t,Y(t))$, we obtain
\begin{align*}
    &\dt \matK_i^< = (\matI \otimes \matQ^\top_{\le i-1}(t_1)) \, [ F(t, Y_i^+(t)) \,]^{\langle i \rangle}  \, \matQ_{\ge i+1}(t_0).
\end{align*}
In an analogous way, the result of the next substep $Y_i^-(t_1)$ is obtained from
\begin{align*}
 &Y_i^-(t) = \matQ_{\le i}(t_1) \, \matS_i(t) \,  \matQ_{\ge i+1}^\top (t_0),  \\
     &\dt \matS_i = - \matQ^\top_{\le i}(t_1)) \, [ F(t, Y_i^-(t)) \,]^{\langle i \rangle}  \, \matQ_{\ge i+1}(t_0).
\end{align*}
These differential equations can be solved numerically by a Runge--Kutta method (of order at least 2 for the second-order splitting integrator). In the important particular case of an autonomous linear ODE 
\[
 \dt A(t) = F(t,A(t)) = L(A(t)), \qquad \text{with linear $L\colon \tensorsize \to \tensorsize$},
\]
the above differential equations are constant-coefficient linear differential equations for $\matK_i^<$ and $\matS_i$, respectively, which can be solved efficiently with a few iterations of a Krylov subspace method for computing the action of the operator exponential~\cite{Hochbruck:1997, sidje-expokit-1998, Hochbruck:2010}.

% \begin{align*}
% \dot \matK_i^< (t) &= (\matI_{n_i} \otimes \matQ_{\le i-1}(\th))^\top \matF^\ivec(t,Y_{i,+}(t))
%  \matQ_{\ge i+1}(t_0)
% \\ & \qquad\hbox{with}\quad \matY_{i,+}^\ivec(t) = (\matI_{n_i} \otimes \matQ_{\le i-1}(\th)) \matK_i^<(t)  \matQ_{\ge i+1}(t_0)^\top
%  \\
% \dot \matS_i(t) &=  -\matQ_{\le i}(\th))^\top \matF^\ivec(t,Y_{i,-}(t))
%  \matQ_{\ge i+1}(t_0)
% \\ & \qquad\hbox{with}\quad \matY_{i,-}^\ivec(t) =  \matQ_{\le i}(\th)) \matS_i(t)  \matQ_{\ge i+1}(t_0)^\top.
% \end{align*}

\subsection{Quantum dynamics in a model potential}
Quantum molecular dynamics is one of the promising applications of the split projector integrator. As a test problem, we use the same setup as considered in~\cite{meyer-henon-2002}:
the time-dependent Schr\"odinger equation with Henon--Heiles potential modeling a coupled oscillator,
\begin{equation}\label{ksl:schrdr}
   i \frac{d\psi}{dt} =  H \psi, \quad \psi(0) = \psi_0,
   \end{equation}
where the Hamiltonian operator $H$ has the form
\begin{equation}\label{eq:hh}
  \def\Hlap{-\frac12 \Delta}
  \def\Hhar{\frac12 \sum_{k=1}^f q^2_k}
  \def\Hanh{\sum_{k=1}^{f-1}\left(q^2_k q_{k+1} - \frac13 q^3_{k+1} \right)}
    H(q_1,\ldots, q_f) = \lefteqn{\overbrace{\phantom{\Hlap + \Hhar}}^{\textrm{harmonic part}}}\Hlap + \underbrace{\Hhar + \overbrace{\lambda \Hanh}^{\textrm{anharmonic part}}}_{\textrm{Henon-Heiles potential}~V(q_1,\ldots,q_f)}
\end{equation}
with $\lambda=0.111803$. 
As an initial condition $\psi_0$, we choose a product of shifted Gaussians,
$$
   \psi_0 = \prod_{i=1}^f \exp\left(-\frac{(q - 2)^2}{2}\right).
$$

The correct discretization of such problems is delicate. A standard approach is to use a Discrete Variable Representation (DVR), specifically, the Sine-DVR scheme from \cite{colbert-sinedvr-1992}.
In addition, since the problem is defined over the whole 
space, appropriate boundary conditions are required. We use complex absorbing potentials (CAP) of the form (see, for example, 
\cite{meyer-book-2009})
  $$
  W(q) = i \eta \sum_{i=1}^f \Big((q_i - q_i^{(r)})_{+}^{b_r} + (q_i - q_i^{(l)})_{-}^{b_l}\Big),
  $$
 where 
 $$
    z_+ = \begin{cases} 
              z, \quad \mbox{if} \quad z \geq 0, \\
              0, \quad \mbox{otherwise},
          \end{cases} \quad \mbox{and} \quad 
    z_- = \begin{cases} 
              z, \quad \mbox{if} \quad z \leq 0, \\
              0, \quad \mbox{otherwise}.
          \end{cases}
 $$
 The parameters $q^{(r)}_i$ and $q^{(l)}_i$ specify the effective boundary of the domain.
 CAP reduces the reflection from the boundary back to the domain, but the system is no longer
 conservative. For the Henon--Heiles example from above we have chosen 
 $$
    \eta = -1, \quad q^{(l)}_i = -6, \quad q^{(r)}_i = 6, \quad b_r = b_l = 3.
 $$

%Our goal is to compute the approximation to the solution of the t

%In order to compute the ground state of $H$, one can consider imaginary time-evolution 
%\[
%???
%\]
%As initial condition, one typically uses shifted Gaussian wave packets since when unshifted these are the exact ground states of harmonic part of $H$.

%The wavefunction should be bounded, which gives the boundary conditions $\psi(\q) \rightarrow 0$, $|\q| \rightarrow \infty$.
%The principal part of the Henon-Heiles operator describes a harmonic oscillator, whose eigenstates are products of a Gaussian by Hermite polynomials and have tensor rank one.
%Therefore, for moderate anharmonicity $\lambda$ the rank one basis of eigenfunctions of the harmonic oscillator is a good choice for the discretization of \eqref{eq:hh}.
%The Galerkin discretization scheme results not only in dense stiffness and mass matrices, but also in a dense matrix that describes the action of the potential $V$.
%The DVR(discrete variable representation) scheme uses a collocation
%of $V$ and $\psi$ on the nodes of the Hermite polynomials and is
%known to provide the same order of accuracy for the eigenproblem as
%the Galerkin method.

We compute the dynamics using the second-order splitting 
integrator where the (linear) local problems for $K_i, \matS_i$ 
are integrated using the Expokit package~\cite{sidje-expokit-1998} with a relative accuracy of $10^{-8}$.

In order to evaluate the accuracy and efficiency 
of our proposed splitting integrator,  we  performed a preliminary comparison with the multi-configuration time-dependent Hartree (MCTDH)
package \cite{mctdh:package}. The MCTDH
method \cite{meyer-book-2009} is the de-facto standard for doing high-dimensional 
quantum molecular dynamics simulations. For the detailed description of MCTDH, we refer to~\cite{meyer-book-2009,
meyer-mctdh-2003,bec-mctdhrev-2000,man-mctdhwav-1992}.

As numerical experiment, we run MCTDH for the 10-dimensional Henon--Heiles problem from above with mode-folding. This can be considered as a first step 
of the hierarchical Tucker format (in this context called the multilayer MCTDH decomposition) with $32$ basis functions in each mode, and the resulting function
was approximated by a 5-dimensional tensor with mode sizes equal to $18$. The final time was $T = 60$.  Our splitting integrator solved the same Henon--Heiles problem but now using the second-order splitting integrator with
a fixed time step $h = 0.01$. Except that we use the TT/MPS manifold for our scheme instead of a Tucker-type manifold as in MCDTH, all other computational parameters are the same.
 
% The experimental setup for our splitting integrator is the same
% as for MCTDH benchmark except that we used a. and d instea  we use another manifold (tensor train) instead
% of the Tucker-type manifold, and we use another integration scheme for the dynamical low-rank approximation which does
% not require any regularization.

\begin{figure}[H]
%\resizebox{13cm}{!}{
\begin{center} \scalebox{0.9}{\input{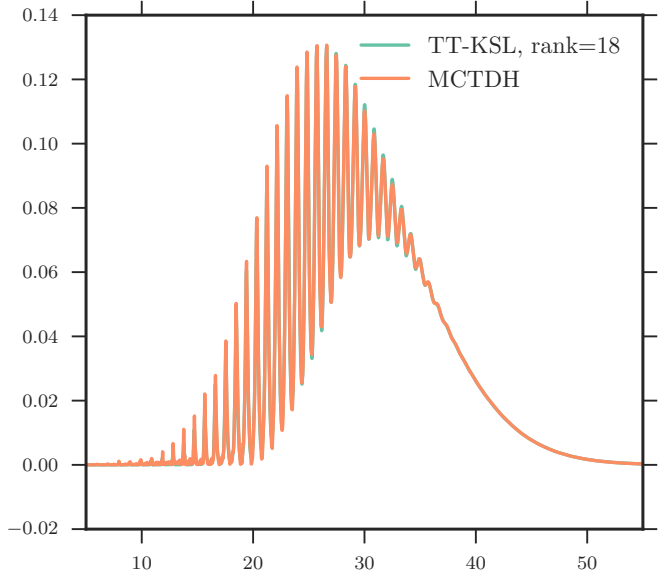}            } \end{center}
%}
%\vspace{-1cm}
\caption{Spectrum computed by the second-order splitting integrator and by the MCTDH package}\label{ksltt:spectra}
\end{figure}

 In Fig.~\ref{ksltt:spectra} we see the vibrational spectrum of  a molecule, which is obtained as follows. After the dynamical low-rank approximation ${\psi}(t)$ is
computed, we evaluate the autocorrelation function 
$
a(t) = \langle \psi(t), \psi(0) \rangle,
$
and compute its Fourier transform $\widehat a(\xi)$. The absolute value
of $\widehat a(\xi)$ gives the information about the energy spectrum of the
operator.  If the dynamics is approximated sufficiently accurately,
the function $\widehat a(\xi)$ is approximated as a sum of delta functions located at 
the eigenvalues of $H$. This method can be considered as a method to approximate
many eigenvalues of $H$ by using only one solution of the dynamical problem, which is not
typical to standard numerical analysis, but often used in chemistry.

We see in~Fig.~\ref{ksltt:spectra} that the computed spectra are very similar, but the MCTDH computation took 
$54\,354$ seconds, whereas the splitting integrator scheme took only $4\,425$ seconds. A detailed comparison of the splitting scheme and MCTDH for quantum molecular dynamics will be presented elsewhere. This will include different benchmark problems and a comparison with the multilayer version of the MCTDH.

\subsection{Approximate matrix inversion}
Optimization on low-rank tensor manifolds is another promising
application of the splitting integrator scheme and can be rather easily incorporated. Consider some iterative process of the form
\begin{equation}\label{ksl:opt1}
Y_{k+1} = Y_k + \Delta_k \quad k = 0, \ldots
\end{equation}
where $\Delta_k$ is the update. In order to obtain approximations $Z_k \in \mathcal{M}$ of $Y_k$ in the TT/MPS format, one typically \emph{retracts} the new iterate back to $\mathcal{M}$,
$$
   Z_{k+1} = P_r( Z_k + \Delta_k),
$$
with $P_r \colon \tensorsize \to \mathcal{M}$ a retraction; see \cite{AMS2008}. A widely used choice for $P_r$ is the quasi-optimal projection computed by TT-SVD~\cite{Oseledets:2011a}. Instead, we propose the cheaper alternative of one step of Algorithm~\ref{al:int} with $A(t_1) - A(t_0) =
\Delta_k$ as $P_r$. In practice, the
intermediate quantities in  Algorithm~\ref{al:int} have to be computed without forming $\Delta_k$ explicitly. This can be done, for example, when $\Delta_k$ is a TT/MPS tensor of low-rank as explained in \S\ref{sec:impl_details}.

An important
example of \eqref{ksl:opt1} is the Newton--Schultz iteration for the approximate matrix
inversion (see, e.g.,~\cite{Hackbusch:2008}),
\begin{equation}\label{ksl:newton}
Y_{k+1} = 2 Y_k - Y_k A Y_k, \quad k = 0, \ldots.
\end{equation}
It is well-known that iteration \eqref{ksl:newton} converges
quadratically provided that $\rho(I - A Y_0) \leq 1$, where
$\rho(\cdot)$ is the spectral radius of the matrix.  The matrix $A$ is
supposed to have low TT/MPS rank when seen as a tensor. This typically arises from a discretization of
a high-dimensional operator on a tensor grid. In our numerical
experiments we have taken the $M$-dimensional Laplace operator with
Dirichlet boundary conditions, discretized by the usual second-order central finite difference on a uniform grid with  $2^d$
points in each mode. 

As a low-rank format, we used the quantized
TT-format (QTT) \cite{osel-2d2d-2010, khor-qtt-2011} which coincides with a $Md$-dimensional TT/MPS format with all dimensions $n_i=2$.  It is known
\cite{khkaz-lap-2012} that in this format the matrix $A$ is represented with
QTT-ranks bounded by $4$. Since $A$ is symmetric positive definite,  as
an initial guess we choose $Y_0 = \alpha I$ with a sufficiently small
$\alpha$.  The splitting integrator is applied with $\Delta_k = Y_k - Y_k A
Y_k$. It requires a certain amount of technical work to implement all
the operations involved in the QTT format, but the final complexity is
linear in the dimension of the tensor (but of course, has high
polynomial complexity with respect to the rank).  
To put the solution
onto the right manifold we artificially add a zero tensor 
%of corresponding rank 
to the initial guess, which has rank $1$, and formally apply the splitting integrator. %This
%brings us to a singular case which would be a disaster for all other
%methods, but the TT-KSL scheme is not influenced by this singularity. 

\begin{figure}[H]
\begin{center}
%\resizebox{10cm}{!}{\input{./ksl-vs-newton-M_2-r_20.pgf}}
\input{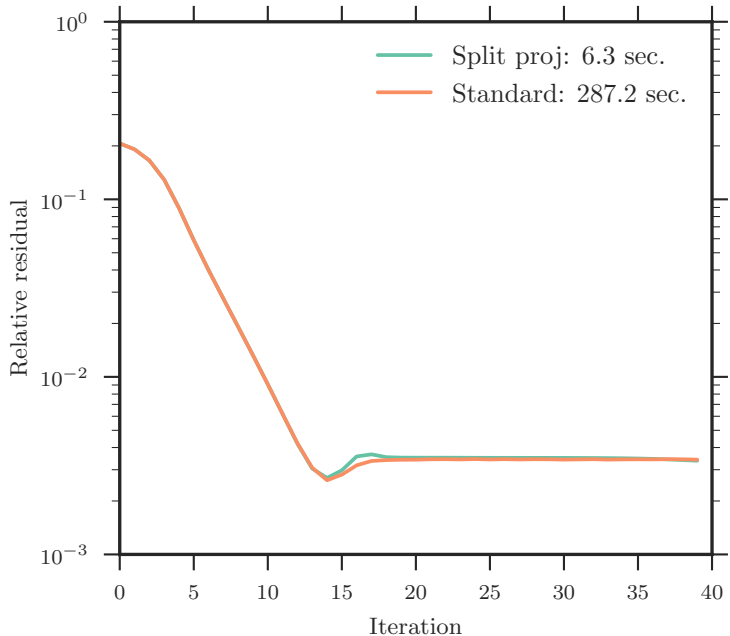}
\caption{Convergence of split projector method and the SVD-based projection method for $D = 2$, $d = 7$, $\alpha = 10^{-2}$, $r = 20$}\label{ksltt:compare}
\end{center}
\end{figure}

\begin{figure}
\input{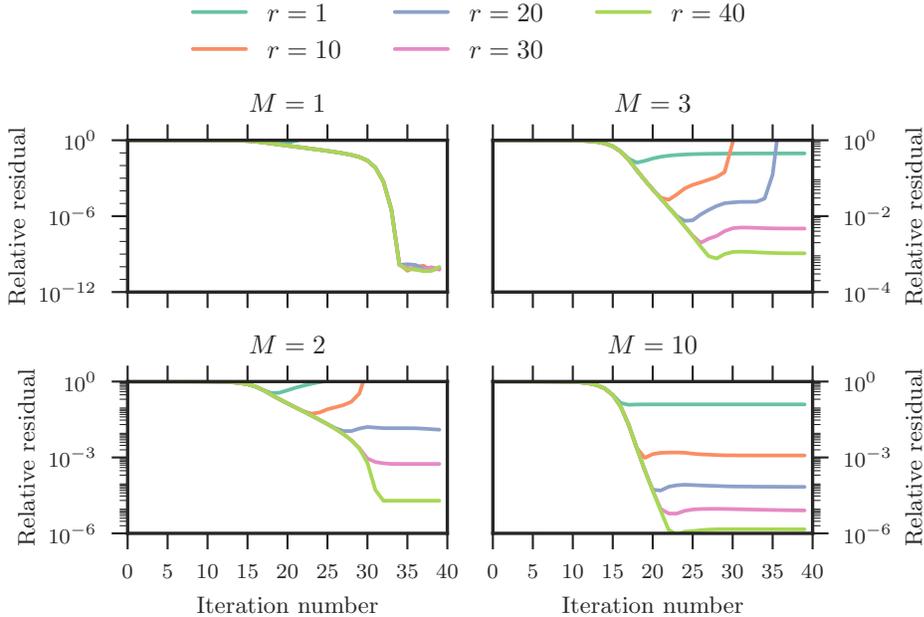}
\caption{The relative residual vs.~iteration number for the approximate
  inversion using TT/MPS rank $r$ of the $M$-dimensional Laplace operator on a uniform grid with $2^7=128$ points in each dimension. 
  Fixed starting guess $Y_0 = \alpha I$ with $\alpha=10^{-6}$.}\label{ksltt:Mrdep}
\end{figure}

\begin{figure}
\input{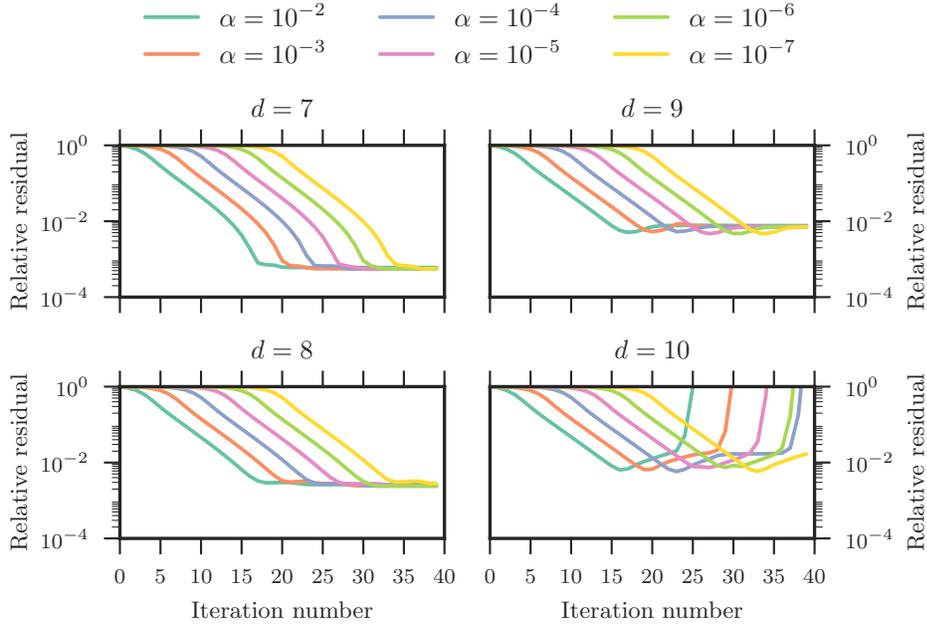}
\caption{The relative residual vs.~iteration number for the approximate
  inversion  using TT/MPS rank $30$ of the $2$-dimensional Laplace operator a uniform grid with $2^d$ points  in each dimension. Starting guesses are $Y_0 = \alpha I$.}\label{ksltt:dalpha}
\end{figure}

\begin{figure}
\begin{center}
\input{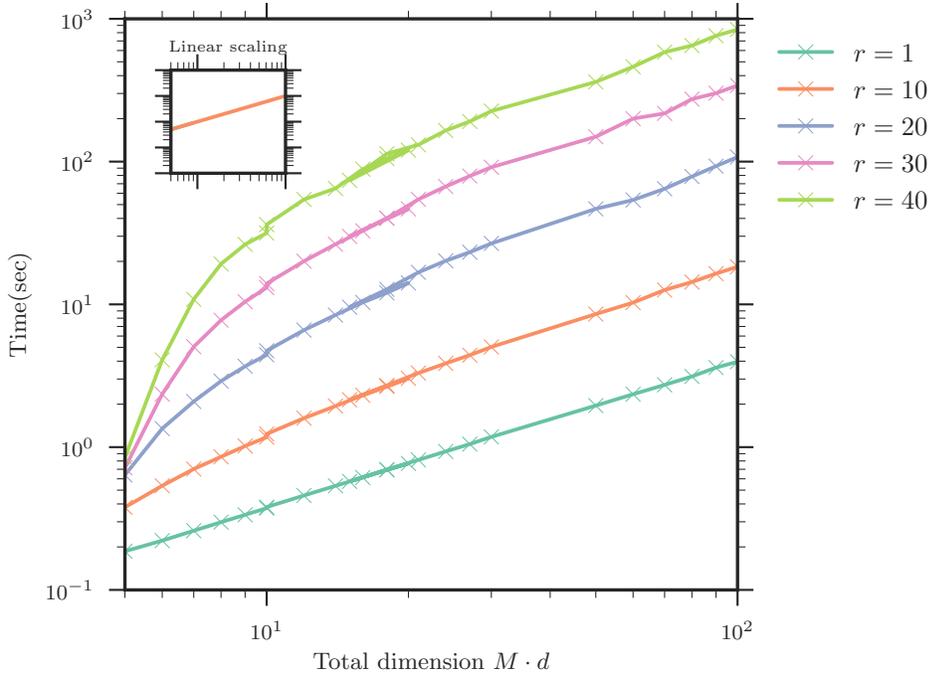}
\caption{Time in log-log scale as a function of the total dimension $M
d$ of the tensor}\label{ksltt:timeschultz}
\end{center}
\end{figure}

As first numerical result, we compare the projector-splitting scheme to the standard approach where after each step of the Newton-Schultz iteration
we project onto a manifold of tensors with bounded TT/MPS ranks $r$ using the TT-SVD,
$$
   Y_{k+1} = P_r(2 Y_k - Y_k A Y_k).
$$
The parameters are set as $M = 2$, $d = 7$, $r = 20$, $\alpha = 10^{-2}$. The convergence of the relative residual $\Vert A
Y_k - I\Vert / \Vert A Y_0 - I \Vert$ in the Frobenius norm for the two methods is presented in Fig.~\ref{ksltt:compare}.
 The splitting method has slightly better accuracy and, more importantly, is significantly faster.
 %The total computation for $40$ iterations
%took $5$ seconds for the split projector whereas standard projection took $553$ seconds.
During the
numerical experiments we observed that the
residual always decreases until the point when the manifold is insufficient to
hold a good approximation to an inverse, and then it either stabilizes
or diverges. The exact explanation of this behavior is out of the scope of the
current paper but could probably be solved using a proper line-search on $\mathcal{M}$ as in~\cite{AMS2008}.  Fig.~\ref{ksltt:Mrdep} shows the convergence
behavior for different $M$ and $r$, with $d$ and $\alpha$ fixed. Fig.~\ref{ksltt:dalpha} shows
the convergence behavior with respect to different $\alpha$ and
$d$. Finally, Fig.~\ref{ksltt:timeschultz} shows that the code has
good scaling with $d$ and $M$.

% Don't forgot to change the first line of conv_from_M.pgf into \pgfpathrectangle{\pgfpointorigin}{\pgfqpoint{8.000000in}{4.500000in}}%
% Do not need now -- the matplotlib code is fixed

% Don't forgot to change the first line of convalphaR30M2.pgf into \pgfpathrectangle{\pgfpointorigin}{\pgfqpoint{8.000000in}{4.650000in}}%

\section{Conclusion}

We have presented and studied a robust and computationally efficient integrator for updating tensors in the tensor train or matrix product state format and for approximately solving tensor differential equations with the approximations retaining the data-sparse tensor train format. Quantum dynamics and tensor optimization appear as promising application areas.

It appears possible to extend this approach to the manifold of hierarchical Tucker tensors of fixed rank \cite{UschmV:2013} and its dynamical approximation \cite{Lubich:2013}. This will be reported elsewhere. In addition, the integrator shares a close resemblance to alternating least squares (ALS) or one-site DMRG (see, e.g., \cite{Holtz:2012,DoOs-dmrg-solve-2011} and for a geometric analysis \cite{Rohwedder:2013}) when the time step goes to infinity. This requires further investigation.

 \section*{Acknowledgement}
 We thank Jutho Haegeman and Frank Verstraete (Gent) for helpful discussions regarding matrix product states and the splitting integrator, and Hans-Dieter Meyer (Heidelberg) for explaining the basic concepts behind quantum molecular dynamics 
 simulations and for his help with the MCTDH package. 
 
We thank the two referees as well as Emil Kieri (Uppsala) and Hanna Walach (T\"ubingen) for pointing out numerous typos in a previous version and for suggesting improvements of the presentation.
 
 The work of C.L. was supported by DFG through SPP 1324 and GRK 1838. The work of I.O. was supported by Russian Science Foundation grant 14-11-00659.

\bibliography{dlrtt,bibtex/molecular,bibtex/misc,bibtex/our,bibtex/tensor,bibtex/dmrg}
\bibliographystyle{siam}

\end{document}